\documentclass[12pt]{amsart}
\usepackage{latexsym, amsmath, amssymb, amsthm}
\usepackage{epsfig}
\usepackage{amscd}
\usepackage{latexsym}
\usepackage{pstricks,pst-node,cite}

\newtheorem{thm}{Theorem}[section]

\newtheorem{cor}[thm]{Corollary}
\newtheorem{exa}[thm]{Example}
\newtheorem{que}[thm]{Question}
\theoremstyle{remark}

\newcommand{\F}{\mathcal{F}}

\newcommand{\D}{\mathcal{D}}

\newcounter{fignum}
\ifx\blackandwhite\undefined
  \newrgbcolor{darkred}{.9 0 0}\newrgbcolor{emgreen}{0 .9 0}
  \newrgbcolor{pup}{.7  0 .9}\newrgbcolor{xxx}{.8 .8 .8}
  \newrgbcolor{lightgray}{.99 .99 0}\newrgbcolor{darkgray}{.8 .8 .3} \else
\fi
\psset{linecolor=blue}

\psset{linewidth=.7pt,dash=3pt 3pt,doublesep=.05,dotsize=1pt 5}
\SpecialCoor

\begin{document}

\author{Dongseok Kim}
\address{Department of Mathematics \\Kyonggi University
\\ Suwon, 443-760 Korea}
\email{dongseok@kgu.ac.kr}

\title[The boundaries of dipole graphs and the complete bipartite graphs $K_{2,n}$]
{The boundaries of dipole graphs and the complete bipartite graphs $K_{2,n}$}

\begin{abstract}
We study the Seifert surfaces of a link by relating the embeddings of graphs by using induced
graphs. As applications, we prove that every link $L$ is the boundary of an oriented surface
which is obtained from a graph embedding of a complete bipartite graph $K_{2,n}$,
where all voltage assignments on the edges of $K_{2,n}$ are $0$.
We also provide an algorithm to construct such a graph diagram of a given link and demonstrate the
algorithm by dealing with the links $4_1^2$ and $5_2$.
\end{abstract}

\maketitle
\section{Introduction}

A \emph{link} $L$ is an embedding of $n$ copies of $\mathbb{S}^1$ into $\mathbb{S}^3$. If a link has only one copy of $\mathbb{S}^1$, the link is called a \emph{knot}.
Throughout the article, we will assume all links are \emph{tame}, which means all links can be deformed in a form of a finite union of line segments.
In the language of graph theory, a knot can be considered as a spatial graph of a cycle graph $C_n$ on $n$ vertices. Two links are \emph{equivalent}
if there is an isotopy between them. In the case of prime knots, this equivalence is the same as the existence of an orientation preserving homeomorphism on $\mathbb{S}^3$,
which sends a knot to the other knot. Although the equivalent class of a link $L$ is called a \emph{link type}, throughout the article, a link
really means the equivalent class of link $L$. Additional terms in the knot theory can be found in~\cite{BZ:knot}.

A \emph{graph} $\Gamma$ is an ordered pair $\Gamma = (V(\Gamma), E(\Gamma))$ comprising a set
$V(\Gamma)$ of vertices together with a set $E(\Gamma)$ of edges. Two graphs $\Gamma_1 = (V(\Gamma_1), E(\Gamma_1))$
and $\Gamma_2 = (V(\Gamma_2), E(\Gamma_2))$ are \emph{equivalent} if there exists a bijective function $\phi :  V(\Gamma_1) \longrightarrow V(\Gamma_2)$
such that $e =\{ u, v \} \in E(\Gamma_1)$ if and only if $\{ \phi(u), \phi(v) \} \in E(\Gamma_2)$.

Although these two subjects are easily considered independent, there are a few branches of graph theory and knot theory which
overlap each other~\cite{DFKLS, Kauffman:jones, Kim:flat, Thistlethwaite, Seifert:def, We}. One of the branches we are interested in is about the Seifert surfaces of the links.
A compact orientable surface $\F$ is called a \emph{Seifert surface} of a link $L$ if the boundary of $\F$ is isotopic to
$L$. The existence of such a surface was first proven by Seifert
using an algorithm on a diagram of $L$. This algorithm was named after him as
\emph{Seifert's algorithm}~\cite{Seifert:def}.
A Seifert surface $\F$ gives rise to a natural signed graph, which is called the \emph{induced graph $\Gamma(\F)$}
by collapsing a disc to a vertex and a twisted band to a signed edge as illustrated in Figure~\ref{graph} $(ii)$.
On the other hand, a graph embedding into an orientable surface naturally produces a Seifert surface by taking a tubular neighborhood of the graph.
These graph embeddings are completely determined by a rotation scheme and a voltage assignment, as explained in Section~\ref{prelim}. However, the objects of having these surfaces
are very different. Knot theorists are concerned about the isotopy classes of surfaces whereas graph theorists only enumerate the homeomorphic classes of graphs in the surfaces.
Besides the sphere, the mapping class groups of orientable surfaces, which are the groups of homeomorphisms on the surface quotient out by the component containing the
identity homeomorphism, are infinite. However, as known as a versatile tool for many different areas,
graphs have been used in many articles related with Seifert surfaces~\cite{JJK:K2n, Kim:flat, Kr, KKL:string, VW}.

In particular, we are interested in banded surfaces and plumbing surfaces. We will relate these surfaces as graph
embeddings of bouquets of circles and dipoles as well as the complete bipartite graphs with integral voltage assignments.
Each surface obtained from the graph embeddings of the bouquets of circles, dipoles and the complete bipartite graphs
is called \emph{bouquet of $n$-circles surface, $n$-dipole surface and the complete bipartite graph surface}, respectively.
If the voltage assignments on all edges are zero, they are said to be \emph{flat}.

A banded surface originally introduced by Kauffman
used to study Seifert pairings, Alexander polynomials of links~\cite{Kauffman}.
The author, Kwon and Lee demonstrated the existence of flat banded surfaces from braid representatives and canonical Seifert surfaces
of a link by using the induced graph of the link~\cite{KKL:string}.

Since the induced graphs of the Seifert surfaces of links are bipartite, all canonical Seifert surfaces of links might be considered as
embeddings of bipartite graphs, where the voltage assignment on the edges of bipartite graphs are $\pm 1$. However, these
graphs are not complete in general.
In a very recent paper, Baader introduced \emph{ribbon diagrams for strongly quasipositive links} \cite{Baader:bipartite}
to show that every $(m, n)$ torus link is a boundary of a surface
which is obtained from the $0$ voltage assignment on all the edges of the complete bipartite graphs $K_{m,n}$ where
the diagram of the complete bipartite graph $K_{m,n}$ is chosen to be in a very special form, as explained as the standard diagram.
For instance, the boundary of the standard diagram of $K_{2,3}$ is the trefoil knot, which is the $(2,3)$ torus knot; yet if we change
a crossing in the standard diagram of $K_{2,3}$, the boundary of this non-standard diagram of $K_{2,3}$ becomes the figure eight knot, which is not
a torus knot; in fact, it is hyperbolic. This phenomenon naturally raises a question whether
every link is a boundary of a complete bipartite graph $K_{n,m}$, as stated in Question~\ref{que1}~\cite{KKL:alternatingsign}.
A weaker version of Question~\ref{que1} was proven immediately without a difficult theory, as in~\cite{JJK:K2n}.
Let us remark that these edges in the graph embeddings are allowed to be linked, but not knotted.
In the present article, we positively answer Question~\ref{que1} that,
for a given link $L$, there exists a graph diagram $D(K_{2,m})$ of a complete bipartite graph $K_{2,m}$
such that link $L$ is a boundary of $D(K_{2,m})$ where all voltage assignments on the edges
of $K_{2,m}$ are $0$. If we flip two discs corresponding to the vertices in one of the bipartition sets whose cardinality is $2$, we obtain
a graph diagram $D(K_{2,m})$ of a complete bipartite graph $K_{2,m}$ whose boundary is the link $L$,
where all voltage assignments on the edges of $K_{2,m}$ are $\pm 1$.

The outline of this paper is as follows. We first provide some preliminary definitions and results
 in Section~\ref{prelim}. In Section\ref{dipole},
we investigate the complete bipartite graph $K_{m,n}$ surface by allowing the bands represented by the edges
of the complete bipartite graph $K_{m,n}$ to be linked but not knotted. In particular, the complete bipartite graph $K_{2,n}$
can be seen as a subdivision of $n$-dipole graphs. We show that every link is a boundary of dipole surface where the signs of the edges are
all $0$(and $\pm 1$). We provide a few examples of such presentations of links.
In Section~\ref{examples}, we provide an algorithm to construct such a graph diagram of a given link and demonstrate the
algorithm by dealing with the links $4_1^2$ and $5_2$.

\section{Preliminaries} \label{prelim}

\begin{figure}
$$
\begin{pspicture}[shift=-2](-.3,-.7)(3.3,4.1) \psline(0,.25)(0,3.75)
\psarc(.25,3.75){0.25}{90}{180}\psline(.25,4)(.5,4)
\pccurve[angleA=0,angleB=180](.5,4)(1,3)
\pccurve[angleA=0,angleB=-110](1,3)(1.227,3.35)
\pccurve[angleA=70,angleB=180](1.3,3.6)(1.5,4) \psline(1.5,4)(2.75,4)
\psarc(2.75,3.75){.25}{0}{90} \psline(3,3.75)(3,.25)
\psarc(2.75,.25){.25}{270}{0} \psline(2.55,0)(2.75,0)
\pccurve[angleA=-50,angleB=180](2.38,.24)(2.55,0)
\pccurve[angleA=-90,angleB=130](2,1)(2.32,.35)
\pccurve[angleA=90,angleB=-90](2,1)(2.5,2) \psline(2.5,2)(2.5,2.75)
\psarc(2.25,2.75){.25}{0}{90}\psline(2.25,3)(2.1,3)
\pccurve[angleA=60,angleB=180](1.82,2.58)(2.1,3)
\pccurve[angleA=90,angleB=-120](1,1)(1.72,2.45)
\pccurve[angleA=135,angleB=-90](1.3,.35)(1,1)
\pccurve[angleA=-45,angleB=180](1.4,.2)(1.75,0)
\psline(2,0)(1.75,0)
\pccurve[angleA=-90,angleB=0](2.5,1)(2,0)
\pccurve[angleA=-60,angleB=90](2.3,1.42)(2.5,1)
\pccurve[angleA=-60,angleB=120](1.6,2.8)(2.2,1.57)
\pccurve[angleA=120,angleB=0](1.6,2.8)(1.5,3)
\pccurve[angleA=0,angleB=180](1,4)(1.5,3)
\pccurve[angleA=60,angleB=-180](.8,3.6)(1,4)
\pccurve[angleA=90,angleB=-120](.5,3)(.7,3.4)
\psline(.5,2.25)(.5,3)
\psarc(.75,2.25){.25}{180}{270}
\psline(.75,2)(.85,2)
\pccurve[angleA=0,angleB=120](.85,2)(1.15,1.7)
\pccurve[angleA=90,angleB=-60](1.5,1)(1.22,1.56)
\pccurve[angleA=0,angleB=-90](1,0)(1.5,1)
\psline(1,0)(.25,0)
\psarc(.25,.25){.25}{180}{270}
\rput[t](.5,1.2){$\alpha$}\rput[t](1.75,1.2){$\beta$}
\rput[t](2.75,2){$\gamma$}\rput[t](1,3.55){$\delta$}
\rput[t](1.5,-.5){$(a)$} \rput[t](-.4,2){$\F_{7_5}$}
\end{pspicture}
\quad\quad  \begin{pspicture}[shift=-2](-.3,-.7)(2.3,4.1) \psline(0,1)(2,1)
\psline(0,1)(1,3) \psline(1,3)(2,1)
\psarc(.5,1){.5}{180}{0}
\psarc(1.5,1){.5}{180}{0}
\pccurve[angleA=70,angleB=110](1,1)(2,1)
\pscircle[fillstyle=solid,fillcolor=lightgray](0,1){.25}
\pscircle[fillstyle=solid,fillcolor=emgreen](1,1){.25}
\pscircle[fillstyle=solid,fillcolor=lightgray](2,1){.25}
\pscircle[fillstyle=solid,fillcolor=emgreen](1,3){.25}
\rput(0,1){$+$} \rput(1,1){$-$} \rput(2,1){$+$} \rput(1,3){$-$}
\rput[t](.5,1.35){$+$}\rput[t](.5,.45){$+$}
\rput[t](1.5,.45){$+$}\rput[t](1.5,.9){$+$}\rput[t](1.5,1.55){$+$}
\rput[t](.4,2.5){$-$}\rput[t](1.6,2.5){$-$}
\rput[t](-.3,1.45){$\alpha$}\rput[t](1,1.65){$\beta$}
\rput[t](2.2,1.45){$\gamma$}\rput[t](1,3.65){$\delta$}
\rput[t](1,-.5){$(b)$} \rput[t](0,3.5){$\Gamma(\F_{7_5})$}
\end{pspicture}
\quad\quad  \begin{pspicture}[shift=-2](-.3,-.7)(2.3,4.1) \psline(0,1)(1,1)
\psline(0,1)(1,3) \psline(1,3)(2,1)
\pscircle[fillstyle=solid,fillcolor=lightgray](0,1){.25}
\pscircle[fillstyle=solid,fillcolor=emgreen](1,1){.25}
\pscircle[fillstyle=solid,fillcolor=lightgray](2,1){.25}
\pscircle[fillstyle=solid,fillcolor=emgreen](1,3){.25}
\rput[t](.5,1.4){$+$}
\rput[t](.4,2.5){$-$}\rput[t](1.6,2.5){$-$}
\rput(0,1){$\alpha$}\rput(1,1){$\beta$}
\rput(2,1){$\gamma$}\rput(1,3){$\delta$}
\rput[t](1,-.5){$(c)$} \rput[t](0,3.5){$T$}
\end{pspicture}
$$
\caption{$(a)$ A knot $7_5$ and its Seifert surface $\F_{7_5}$ whose discs are named
$\alpha, \beta, \gamma, \delta$, $(b)$ its corresponding signed induced graph $\Gamma(\F_{7_5})$ and
$(c)$ a spanning tree $T$ of $\Gamma(\F_{7_5})$.} \label{graph}
\end{figure}
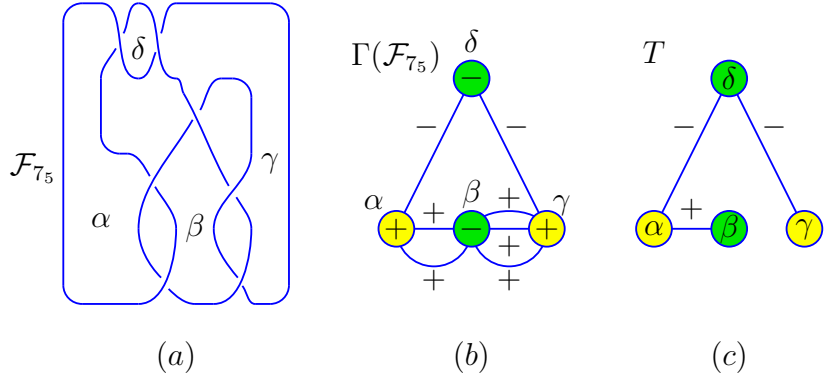

A Seifert surface $\F_L$ of an oriented link $L$ is produced by applying Seifert's
algorithm to a link diagram $D(L)$ as shown in Figure~\ref{graph} $(i)$; it
is called a \emph{canonical Seifert surface}.
A canonical Seifert surface $\F$ gives rise to a natural signed graph, which is referred to as the \emph{induced graph $\Gamma(\F)$}
by collapsing a disc to a vertex and a twisted band to a signed edge, as illustrated in Figure~\ref{graph} $(i)$.
These processes can also be performed on arbitrary Seifert surfaces.
Since link $L$ is tame and its Seifert surface $\F_L$ is compact, the induced graph $\Gamma(\F_L)$ is finite.
By separating the discs by local orientation as indicated $\pm$ on each vertices in Figure~\ref{graph} $(ii)$, the induced graph $\Gamma(\F_L)$ can be considered as a bipartite graph.
If the Seifert surface is connected, then its induced graph is also connected.

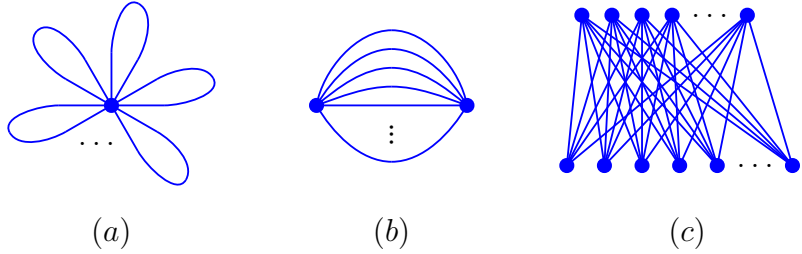
\begin{figure}
$$
\begin{pspicture}[shift=-1.2](-1.6,-2.1)(1.6,1.7)
\pscircle[fillcolor=lightgray, fillstyle=solid, linewidth=3pt](0,0){.1}
\psline(.7;0)(.7;180) \psline(.7;30)(.7;210) \psline(.7;60)(0;240)
\psline(.7;90)(0;270) \psline(.7;120)(.7;300) \psline(.7;150)(.7;330)
\pccurve[angleA=0,angleB=-60](.7;0)(1.4;15)
\pccurve[angleA=30,angleB=120](.7;30)(1.4;15)
\pccurve[angleA=60,angleB=-30](.7;60)(1.4;75)
\pccurve[angleA=90,angleB=150](.7;90)(1.4;75)
\pccurve[angleA=120,angleB=45](.7;120)(1.4;135)
\pccurve[angleA=150,angleB=-135](.7;150)(1.4;135)
\pccurve[angleA=180,angleB=120](.7;180)(1.4;195)
\pccurve[angleA=210,angleB=-60](.7;210)(1.4;195)
\pccurve[angleA=300,angleB=-120](.7;300)(1.4;315)
\pccurve[angleA=330,angleB=60](.7;330)(1.4;315)
\rput(-.2,-.5){{$\ldots$}}
\rput(0,-1.7){{$(a)$}}
\end{pspicture}\quad\quad
\begin{pspicture}[shift=-1.2](-.3,-2.1)(2.3,1.7)
\pscircle[fillcolor=lightgray, fillstyle=solid, linewidth=3pt](0,0){.1}
\pscircle[fillcolor=lightgray, fillstyle=solid, linewidth=3pt](2,0){.1}
\pccurve[angleA=80,angleB=180](0,0)(1,1) \pccurve[angleA=0,angleB=100](1,1)(2,0)
\pccurve[angleA=60,angleB=180](0,0)(1,.75) \pccurve[angleA=0,angleB=120](1,.75)(2,0)
\pccurve[angleA=40,angleB=180](0,0)(1,.5) \pccurve[angleA=0,angleB=140](1,.5)(2,0)
\pccurve[angleA=20,angleB=180](0,0)(1,.25) \pccurve[angleA=0,angleB=160](1,.25)(2,0)
\psline(0,0)(2,0)
\rput(1,-.3){{$\cdot$}} \rput(1,-.4){{$\cdot$}} \rput(1,-.5){{$\cdot$}}
\pccurve[angleA=-60,angleB=180](0,0)(1,-.75) \pccurve[angleA=0,angleB=-120](1,-.75)(2,0)
\rput(1,-1.7){{$(b)$}}
\end{pspicture} \quad\quad
\begin{pspicture}[shift=-1.2](-.2,-2.1)(3.4,1.7)
\pscircle[fillcolor=lightgray, fillstyle=solid, linewidth=3pt](0,-.8){.1}
\pscircle[fillcolor=lightgray, fillstyle=solid, linewidth=3pt](.5,-.8){.1}
\pscircle[fillcolor=lightgray, fillstyle=solid, linewidth=3pt](1,-.8){.1}
\pscircle[fillcolor=lightgray, fillstyle=solid, linewidth=3pt](1.5,-.8){.1}
\pscircle[fillcolor=lightgray, fillstyle=solid, linewidth=3pt](2,-.8){.1}
\pscircle[fillcolor=lightgray, fillstyle=solid, linewidth=3pt](3,-.8){.1}
\pscircle[fillcolor=lightgray, fillstyle=solid, linewidth=3pt](.2,1.2){.1}
\pscircle[fillcolor=lightgray, fillstyle=solid, linewidth=3pt](.6,1.2){.1}
\pscircle[fillcolor=lightgray, fillstyle=solid, linewidth=3pt](1,1.2){.1}
\pscircle[fillcolor=lightgray, fillstyle=solid, linewidth=3pt](1.4,1.2){.1}
\pscircle[fillcolor=lightgray, fillstyle=solid, linewidth=3pt](2.4,1.2){.1}
\psline(0,-.8)(.2,1.2)
\psline(0,-.8)(.6,1.2)
\psline(0,-.8)(1,1.2)
\psline(0,-.8)(1.4,1.2)
\psline(0,-.8)(2.4,1.2)
\psline(.5,-.8)(.2,1.2)
\psline(.5,-.8)(.6,1.2)
\psline(.5,-.8)(1,1.2)
\psline(.5,-.8)(1.4,1.2)
\psline(.5,-.8)(2.4,1.2)
\psline(1,-.8)(.2,1.2)
\psline(1,-.8)(.6,1.2)
\psline(1,-.8)(1,1.2)
\psline(1,-.8)(1.4,1.2)
\psline(1,-.8)(2.4,1.2)
\psline(1.5,-.8)(.2,1.2)
\psline(1.5,-.8)(.6,1.2)
\psline(1.5,-.8)(1,1.2)
\psline(1.5,-.8)(1.4,1.2)
\psline(1.5,-.8)(2.4,1.2)
\psline(2,-.8)(.2,1.2)
\psline(2,-.8)(.6,1.2)
\psline(2,-.8)(1,1.2)
\psline(2,-.8)(1.4,1.2)
\psline(2,-.8)(2.4,1.2)
\psline(3,-.8)(.2,1.2)
\psline(3,-.8)(.6,1.2)
\psline(3,-.8)(1,1.2)
\psline(3,-.8)(1.4,1.2)
\psline(3,-.8)(2.4,1.2)
\rput(1.9,1.2){{$\ldots$}} \rput(2.5,-.8){{$\ldots$}}
\rput(1.6,-1.7){{$(c)$}}
\end{pspicture}
$$
\caption{$(a)$ A bouquet $n$-circles graph $B_n$, $(b)$ an $n$-dipole graph $D_n$
and $(c)$ the complete bipartite graph $K_{m,n}$.} \label{bouquetfig1}
\end{figure}

It is fairly easy to see that the number of Seifert
circles (half twisted bands), denoted by $s(\F_L)(c(\F_L))$, is the
cardinality of the vertex set, $V(\Gamma(\F_L))$ (edge set $E(\Gamma(\F_L))$, respectively).
A spanning tree $T$ of $\Gamma(\F_L)$ is depicted in Figure~\ref{graph} $(iii)$.
Recall that the number of edges of a spanning tree of a connected graph with $n$ vertices is $n-1$.
One can see that the length of the path joining both ends of an arbitrary vertex $e \in \Gamma(\F_{L})$ is odd.

Next, let us provide the definitions of graphs which are used in the article. A \emph{bouquet of $n$-circles}, denoted by $B_n$,
is a graph with a single vertex and $n$ self loops as illustrated in Figure~\ref{bouquetfig1} $(a)$.
These bouquets of circles are fundamental building blocks in topological graph theory
because any connected graph can be reduced to a bouquet of circles by contracting a
spanning tree to a point as all coedges become generators of the fundamental group of the graph.
An \emph{$n$-dipole}, denoted by $D_n$
is a graph with two vertices and $n$ edges joining these two vertices as depicted in Figure~\ref{bouquetfig1} $(b)$.
A complete bipartite graph, $\Gamma = (V_1 \amalg V_2, E)$, is a bipartite graph such that
for any two vertices, $v_1 \in V_1$ and $v_2 \in V_2$, $\{ v_1, v_2 \}$ is an edge in $\Gamma$.
The complete bipartite graph with partitions of size $|V_1|=m$ and $|V_2|=n$, is denoted $K_{m,n}$ as shown in Figure~\ref{bouquetfig1} $(c)$.

\begin{figure}
$$
\begin{pspicture}[shift=-1.2](-2.5,-2)(2.5,2.4)
\psline(0,0)(2;90) \psline(0,0)(2;-30) \psline(0,0)(2;-150)
\psline(2;-30)(2;90) \psline(2;-150)(2;90) \psline(2;-30)(2;-150)
\pscircle[fillcolor=lightgray, fillstyle=solid, linewidth=1pt](0,0){.4}
\pscircle[fillcolor=lightgray, fillstyle=solid, linewidth=1pt](2;90){.4}
\pscircle[fillcolor=lightgray, fillstyle=solid, linewidth=1pt](2;-30){.4}
\pscircle[fillcolor=lightgray, fillstyle=solid, linewidth=1pt](2;-150){.4}
\rput(0,0){{$\alpha$}} \rput(2;90){{$\beta$}} \rput(2;-150){{$\gamma$}} \rput(2;-30){{$\delta$}}
\rput(0,-.6){{$\beta\gamma\delta$}} \rput(.75,2.25){{$\alpha\delta\gamma$}} \rput(2.35,-.55){{$\alpha\gamma\beta$}} \rput(-2.35,-.55){{$\alpha\beta\delta$}}
\rput(.9,-.2){{$1$}} \rput(-.9,-.2){{$1$}} \rput(.8,1.1){{$1$}}
\rput(-.8,1.1){{$1$}} \rput(.2,1){{$0$}} \rput(0,-1.2){{$0$}}
\rput(0,-1.8){{$(a)$}}
\end{pspicture}
\quad
\begin{pspicture}[shift=-1.2](-2.5,-2)(2.5,2.4)
\psline(.095,0)(.095,2) \psline(-.095,0)(-.095,2)
\psline(.05,.0732)(0.7428,-.3268) \psline(1.0892,-.5268)(1.782,-0.9268)
\psline(-.05,-.0732)(.6428,-.4732) \psline(.9892,-.6732)(1.682,-1.0732)
\pccurve[angleA=-30,angleB=150](0.7428,-.3268)(.9892,-.6732)
\pccurve[angleA=-30,angleB=-160](.6428,-.4732)(.826,-.505)
\pccurve[angleA=10,angleB=150](.906,-.495)(1.0892,-.5268)
\psline(-.05,.0732)(-.7428,-.3268) \psline(-1.0892,-.5268)(-1.782,-0.9268)
\psline(.05,-.0732)(-.6428,-.4732) \psline(-.9892,-.6732)(-1.682,-1.0732)
\pccurve[angleA=-150,angleB=30](-0.7428,-.3268)(-.9892,-.6732)
\pccurve[angleA=-150,angleB=-20](-.6428,-.4732)(-.826,-.505)
\pccurve[angleA=170,angleB=30](-.906,-.495)(-1.0892,-.5268)
\psline(-1.732,-.905)(1.732,-0.905)
\psline(-1.732,-1.095)(1.732,-1.095)
\psline(.0732,2.05)(.8526,.7) \psline(1.0258,.4)(1.8052,-.95)
\psline(-.0732,1.95)(.7062,.6) \psline(.8794,.3)(1.6588,-1.05)
\pccurve[angleA=-60,angleB=120](.8526,.7)(.8794,.3)
\pccurve[angleA=-60,angleB=160](.7062,.6)(.826,.505)
\pccurve[angleA=-20,angleB=120](.906,.495)(1.0258,.4)
\psline(-.0732,2.05)(-.8526,.7) \psline(-1.0258,.4)(-1.8052,-.95)
\psline(.0732,1.95)(-.7062,.6) \psline(-.8794,.3)(-1.6588,-1.05)
\pccurve[angleA=-120,angleB=60](-.8526,.7)(-.8794,.3)
\pccurve[angleA=-120,angleB=20](-.7062,.6)(-.826,.505)
\pccurve[angleA=160,angleB=60](-.906,.495)(-1.0258,.4)
\pscircle[fillcolor=lightgray, fillstyle=solid, linewidth=1pt](0,0){.4}
\pscircle[fillcolor=lightgray, fillstyle=solid, linewidth=1pt](2;90){.4}
\pscircle[fillcolor=emgreen, fillstyle=solid, linewidth=1pt](2;-30){.4}
\pscircle[fillcolor=emgreen, fillstyle=solid, linewidth=1pt](2;-150){.4}
\rput(0,0){{$\alpha$}} \rput(2;90){{$\beta$}} \rput(2;-150){{$\gamma$}} \rput(2;-30){{$\delta$}}
\rput(0,-1.8){{$(b)$}}
\end{pspicture}
$$
$$
\begin{pspicture}[shift=-1.2](-2.4,-3.3)(3.4,2.8)
\pccurve[doubleline=true, angleA=45,angleB=150](0,0)(2,-.2)
\pccurve[doubleline=true, angleA=-30,angleB=45](2,-.2)(2.8,-1.5)
\pccurve[doubleline=true, angleA=-135,angleB=-45](2.8,-1.5)(2,-1.2)
\pccurve[doubleline=true, angleA=-90,angleB=0](0,0)(-1.5,-2)
\pccurve[doubleline=true, angleA=180,angleB=-135](-1.5,-2)(-2.1,-1)
\psline[doubleline=true](2;90)(0,0)
\psline[doubleline=true](2;90)(2;-30)
\psline[doubleline=true](2;90)(2;-150)
\psline[doubleline=true](2;-150)(2;-30)
\pscircle[fillcolor=lightgray, fillstyle=solid, linewidth=1pt](0,0){.4}
\pscircle[fillcolor=lightgray, fillstyle=solid, linewidth=1pt](2;90){.4}
\pscircle[fillcolor=lightgray, fillstyle=solid, linewidth=1pt](2;-30){.4}
\pscircle[fillcolor=lightgray, fillstyle=solid, linewidth=1pt](2;-150){.4}
\rput(0,0){{$\alpha$}} \rput(2;90){{$\beta$}} \rput(2;-150){{$\gamma$}} \rput(2;-30){{$\delta$}}
\rput(0,-3){{$(c)$}}
\end{pspicture}
\quad
\begin{pspicture}[shift=-1.2](-2.5,-3.3)(2.5,2.8)
\psline[doubleline=true](0,0)(0,2.414)
\psline[doubleline=true](0,0)(1.707,1.707)
\psline[doubleline=true](0,0)(0,-1)
\pccurve[doubleline=true, angleA=-90,angleB=45](0,-1)(-1.707,-1.707)
\pccurve[doubleline=true, angleA=90,angleB=180](1.5;-30)(2.414,0)
\pccurve[doubleline=true, angleA=-150,angleB=90](1.5;-30)(0,-2.414)
\pccurve[doubleline=true, angleA=-30,angleB=135](1.5;-30)(1.707,-1.707)
\pccurve[doubleline=true, angleA=-120,angleB=0](1.5;-150)(-2.414,0)
\pccurve[doubleline=true, angleA=135,angleB=-45](1.5;-150)(-1.707,1.707)
\psline[doubleline=true](1.5;-150)(1.5;90)
\pscircle[fillcolor=lightgray, fillstyle=solid, linewidth=1pt](0,0){.4}
\pscircle[fillcolor=lightgray, fillstyle=solid, linewidth=1pt](1.5;90){.4}
\pscircle[fillcolor=lightgray, fillstyle=solid, linewidth=1pt](1.5;-30){.4}
\pscircle[fillcolor=lightgray, fillstyle=solid, linewidth=1pt](1.5;-150){.4}
\rput(0,0){{$\alpha$}} \rput(1.5;90){{$\beta$}} \rput(1.5;-150){{$\gamma$}} \rput(1.5;-30){{$\delta$}}
\psline[linecolor=pup, linewidth=1.4pt](1,2.414)(-1,2.414)
\psline[linecolor=pup, arrowscale=2.8]{->}(.4,2.414)(.6,2.414)
\psline[linecolor=emgreen, linewidth=1.4pt](1,2.414)(2.414,1)
\psline[linecolor=emgreen, arrowscale=2.8]{->}(2.0605,1.3535)(2.2605,1.1535)
\psline[linecolor=pup, linewidth=1.4pt](2.414,1)(2.414,-1)
\psline[linecolor=pup, arrowscale=2.8]{->}(2.414,.4)(2.414,.6)
\psline[linecolor=emgreen, linewidth=1.4pt](1,-2.414)(2.414,-1)
\psline[linecolor=emgreen, arrowscale=2.8]{->}(2.0605,-1.3535)(2.2605,-1.1535)
\psline[linecolor=darkred, linewidth=1.4pt](-1,-2.414)(1,-2.414)
\psline[linecolor=darkred, arrowscale=2.8]{->}(-.4,-2.414)(-.6,-2.414)
\psline[linecolor=darkgray, linewidth=1.4pt](-1,-2.414)(-2.414,-1)
\psline[linecolor=darkgray, arrowscale=2.8]{->}(-2.0605,-1.3535)(-2.2605,-1.1535)
\psline[linecolor=darkred, linewidth=1.4pt](-2.414,-1)(-2.414,1)
\psline[linecolor=darkred, arrowscale=2.8]{->}(-2.414,-.4)(-2.414,-.6)
\psline[linecolor=darkgray, linewidth=1.4pt](-1,2.414)(-2.414,1)
\psline[linecolor=darkgray, arrowscale=2.8]{->}(-2.0605,1.3535)(-2.2605,1.1535)
\rput(0,-3){{$(d)$}}
\end{pspicture}
$$
\caption{$(a)$ The complete bipartite graph $K_4$ with a rotation scheme and $\mathbb{Z}_2$ voltage assignment,
$(b)$ a band decomposition of $K_4$ surface with respect to the rotation scheme and $\mathbb{Z}_2$ voltage assignment in $(a)$,
$(c)$ a flat $K_4$ surface obtained from $(b)$ by flipping discs $\gamma$ and $\delta$
and $(d)$ an embedding of the complete bipartite graph $K_4$ into genus $2$ torus which is obtained from $(c)$ by attaching the two
discs along the boundary of a flat $K_4$ surface.} \label{rotation}
\end{figure}
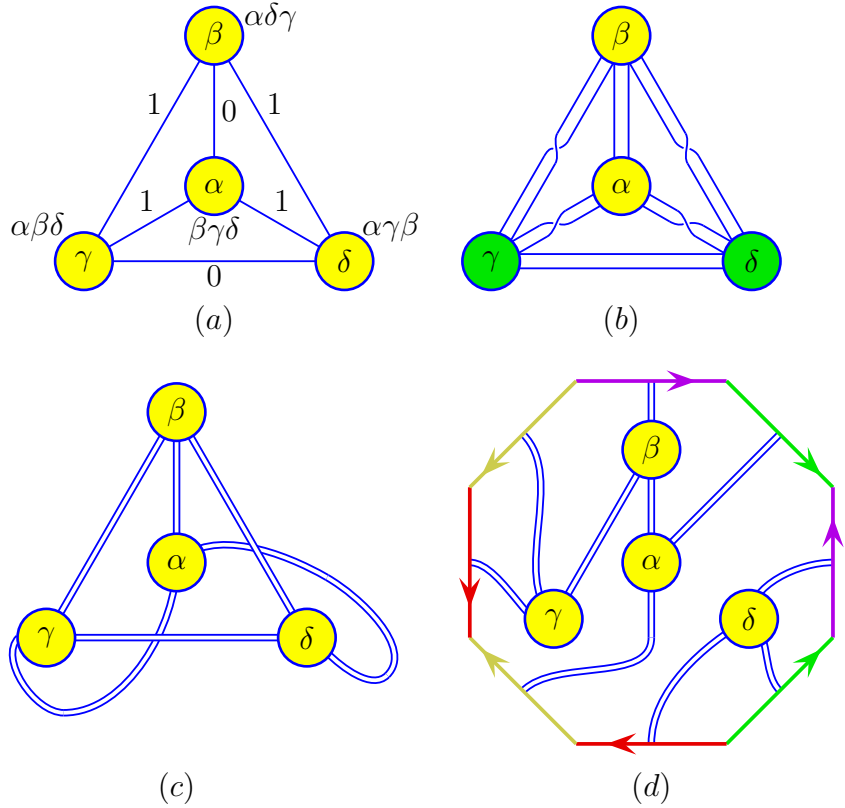

Let $\Gamma=(V(\Gamma), E(\Gamma))$ be a graph.
A \emph{rotation scheme} is a set of cyclic orders of edges which are adjacent to $v$ for all vertex $v\in V(\Gamma)$,
A \emph{voltage assignment} is a map from $E(\Gamma)$ to $\mathbb{Z}_2$. These information are often called a \emph{map}.
Next, we will explain how the rotation scheme and embeddings of a graph $\Gamma$ determine the embedding of $\Gamma$ into an orientable surface.
The main idea of the construction of a surface is a band decomposition, as follows:
each vertex $v \in E(\Gamma)$ is replaced by a disc as $0$-band,
$1$-bands representing an edge are glued in order to respect the rotation scheme(cyclic order on each vertices) and the voltage assignment
(if the value on an edge is $0$, the band is flat, otherwise, the band is half twisted)
and $2$-bands are capping off the boundaries~\cite{GT1}. Figure~\ref{rotation} provides an example
where $\Gamma$ is the complete graph $K_4$. Let us remark that the boundary of the flat $K_4$ surface in Figure~\ref{rotation} $(c)$
is the hope link.

As described earlier, a Seifert surface produces the induced graph. Conversely, a graph produces
a surface which may not be orientable. In order to have the orientability of a resulting surface,
we will assume that the voltage assignment is zero for all edges. In graph embeddings, graph theorists are mainly concerned about
the enumerations of non equivalent embeddings of a graph into a surface. In contrast, knot theorists are concerned about the isotopy classes instead of homeomorphism classes of surfaces.
In fact, there exist non equivalent links which are boundaries of homeomorphic surfaces obtained from the two cell embeddings of the same graph.
Figure~\ref{k23torus} and Figure~\ref{k23figure8}, which are two cell embeddings of the graph $K_{2,3}$, have different boundaries, the trefoil and the figure eight knot.
The relation between graph embeddings and the Seifert surface are studied extensively with topological graph theorists~\cite{BKKL}.

\begin{figure}
$$
\begin{pspicture}[shift=-1.3](-1.4,-1.7)(1.4,2.1)
\pscircle[fillcolor=lightgray, fillstyle=solid, linewidth=3pt](.4,0){.1}
\pscircle[linewidth=3pt](1.2,0){.1}
\pscircle[linewidth=3pt](-1.2,0){.1}
\pscircle[linewidth=3pt](-.4,0){.1}
\pccurve[linewidth=1.2pt, angleA=60,angleB=-120](-1.2,0)(-.84,.54)
\pccurve[linewidth=1.2pt, angleA=60,angleB=0](-.76,.66)(-.75,1.45)
\pccurve[linewidth=1.2pt, angleA=180,angleB=120](-.85,1.45)(-.8,.6)
\pccurve[linewidth=1.2pt, angleA=-60,angleB=150](-.8,.6)(-.05,.33)
\pccurve[linewidth=1.2pt, angleA=-30,angleB=120](.05,.27)(.4,0)
\pccurve[linewidth=1.2pt, angleA=60,angleB=-150](-.4,0)(0,.3)
\pccurve[linewidth=1.2pt, angleA=30,angleB=-120](0,.3)(.36,.54)
\pccurve[linewidth=1.2pt, angleA=60,angleB=-60](.44,.66)(.42,1.4)
\pccurve[linewidth=1.2pt, angleA=120,angleB=0](.37,1.5)(0,2)
\pccurve[linewidth=1.2pt, angleA=180,angleB=90](0,2)(-.8,1.45)
\pccurve[linewidth=1.2pt, angleA=-90,angleB=170](-.8,1.45)(-.65,1.125)
\pccurve[linewidth=1.2pt, angleA=-10,angleB=-170](-.55,1.1)(.45,1.1)
\pccurve[linewidth=1.2pt, angleA=10,angleB=-90](.57,1.13)(1.1,1.5)
\pccurve[linewidth=1.2pt, angleA=90,angleB=0](1.1,1.5)(.9,1.8)
\pccurve[linewidth=1.2pt, angleA=180,angleB=50](.9,1.8)(.4,1.45)
\pccurve[linewidth=1.2pt, angleA=-130,angleB=70](.4,1.45)(.21,1.13)
\pccurve[linewidth=1.2pt, angleA=-110,angleB=135](.155,1.0)(.4,.6)
\pccurve[linewidth=1.2pt, angleA=-45,angleB=120](.4,.6)(1.2,0)
\rput(0,-1.5){{$(a)$}}
\end{pspicture} \quad
\begin{pspicture}[shift=-1.3](-2.2,-1.7)(2.2,2.1)
\pccurve[doubleline=true, angleA=90,angleB=-120](-1.2,0)(-.7,.8)
\pccurve[doubleline=true, angleA=60,angleB=0](-.7,.8)(-.7,1.45)
\pccurve[doubleline=true, angleA=180,angleB=0](-.7,1.45)(-.9,1.45)
\pccurve[doubleline=true, angleA=180,angleB=120](-.9,1.45)(-.8,.6)
\pccurve[doubleline=true, angleA=-60,angleB=150](-.8,.6)(0,.3)
\pccurve[doubleline=true, angleA=-30,angleB=90](0,.3)(.4,0)
\pccurve[doubleline=true, angleA=90,angleB=-150](-.4,0)(0,.3)
\pccurve[doubleline=true, angleA=30,angleB=-120](0,.3)(.5,.7)
\pccurve[doubleline=true, angleA=60,angleB=-90](.5,.7)(.6,1.3)
\pccurve[doubleline=true, angleA=90,angleB=-60](.6,1.3)(.4,1.8)
\pccurve[doubleline=true, angleA=120,angleB=0](.4,1.8)(0,2)
\pccurve[doubleline=true, angleA=180,angleB=90](0,2)(-.8,1.45)
\pccurve[doubleline=true, angleA=-90,angleB=170](-.8,1.45)(-.6,1.1)
\pccurve[doubleline=true, angleA=-10,angleB=180](-.6,1.1)(0,1.05)
\pccurve[doubleline=true, angleA=0,angleB=-160](0,1.05)(.8,1.2)
\pccurve[doubleline=true, angleA=20,angleB=-90](.8,1.2)(1.1,1.5)
\pccurve[doubleline=true, angleA=90,angleB=0](1.1,1.5)(.9,1.8)
\pccurve[doubleline=true, angleA=180,angleB=50](.9,1.8)(.4,1.45)
\pccurve[doubleline=true, angleA=-130,angleB=70](.4,1.45)(.2,1)
\pccurve[doubleline=true, angleA=-110,angleB=135](.2,1)(.4,.6)
\pccurve[doubleline=true, angleA=-45,angleB=90](.4,.6)(1.2,0)
\pccurve[doubleline=true, angleA=60,angleB=0](-.7,.8)(-.7,1.45)
\pccurve[doubleline=true, angleA=0,angleB=-160](0,1.05)(.8,1.2)
\pccurve[doubleline=true, angleA=60,angleB=-90](.5,.7)(.6,1.3)
\psframe[linecolor=white, fillstyle=solid, fillcolor=lightgray](-1.6,-1)(1.6,.1)
\psline(-1.215,.1)(-1.6,.1)(-1.6,-1)(1.6,-1)(1.6,.1)(1.217,.1)
\psline(-1.156,.1)(-.395,.1) \psline(-.317,.1)(.318,.1) \psline(1.148,.1)(.4,.1)
\rput(0,-1.5){{$(b)$}}
\rput(0,-.5){{$\mathcal{D}$}}
\end{pspicture}
\quad
\begin{pspicture}[shift=-1.3](-1.4,-1.7)(1.4,2.1)
\pscircle[fillcolor=lightgray, fillstyle=solid, linewidth=3pt](0,-.4){.1}
\pccurve[linewidth=1.2pt, angleA=-120,angleB=-135](-1.2,0)(0,-.4)
\pccurve[linewidth=1.2pt, angleA=60,angleB=-120](-1.2,0)(-.84,.54)
\pccurve[linewidth=1.2pt, angleA=60,angleB=0](-.76,.66)(-.75,1.45)
\pccurve[linewidth=1.2pt, angleA=180,angleB=120](-.85,1.45)(-.8,.6)
\pccurve[linewidth=1.2pt, angleA=-60,angleB=150](-.8,.6)(-.05,.33)
\pccurve[linewidth=1.2pt, angleA=-30,angleB=120](.05,.27)(.4,0)
\pccurve[linewidth=1.2pt, angleA=-60,angleB=45](.4,0)(0,-.4)
\pccurve[linewidth=1.2pt, angleA=-120,angleB=135](-.4,0)(0,-.4)
\pccurve[linewidth=1.2pt, angleA=60,angleB=-150](-.4,0)(0,.3)
\pccurve[linewidth=1.2pt, angleA=30,angleB=-120](0,.3)(.36,.54)
\pccurve[linewidth=1.2pt, angleA=60,angleB=-60](.44,.66)(.42,1.4)
\pccurve[linewidth=1.2pt, angleA=120,angleB=0](.37,1.5)(0,2)
\pccurve[linewidth=1.2pt, angleA=180,angleB=90](0,2)(-.8,1.45)
\pccurve[linewidth=1.2pt, angleA=-90,angleB=170](-.8,1.45)(-.65,1.125)
\pccurve[linewidth=1.2pt, angleA=-10,angleB=-170](-.55,1.1)(.45,1.1)
\pccurve[linewidth=1.2pt, angleA=10,angleB=-90](.57,1.13)(1.1,1.5)
\pccurve[linewidth=1.2pt, angleA=90,angleB=0](1.1,1.5)(.9,1.8)
\pccurve[linewidth=1.2pt, angleA=180,angleB=50](.9,1.8)(.4,1.45)
\pccurve[linewidth=1.2pt, angleA=-130,angleB=70](.4,1.45)(.21,1.13)
\pccurve[linewidth=1.2pt, angleA=-110,angleB=135](.155,1.0)(.4,.6)
\pccurve[linewidth=1.2pt, angleA=-45,angleB=120](.4,.6)(1.2,0)
\pccurve[linewidth=1.2pt, angleA=-60,angleB=-45](1.2,0)(0,-.4)
\rput(0,-1.5){{$(c)$}}
\end{pspicture} \quad
\begin{pspicture}[shift=-1.3](-2.2,-1.7)(2.2,2.1)
\pccurve[doubleline=true, angleA=-135,angleB=-90](0,-.4)(-1.2,0)
\pccurve[doubleline=true, angleA=90,angleB=-120](-1.2,0)(-.7,.8)
\pccurve[doubleline=true, angleA=60,angleB=0](-.7,.8)(-.7,1.45)
\pccurve[doubleline=true, angleA=180,angleB=0](-.7,1.45)(-.9,1.45)
\pccurve[doubleline=true, angleA=180,angleB=120](-.9,1.45)(-.8,.6)
\pccurve[doubleline=true, angleA=-60,angleB=150](-.8,.6)(0,.3)
\pccurve[doubleline=true, angleA=-30,angleB=90](0,.3)(.4,0)
\pccurve[doubleline=true, angleA=-90,angleB=45](.4,0)(0,-.4)
\pccurve[doubleline=true, angleA=135,angleB=-90](0,-.4)(-.4,0)
\pccurve[doubleline=true, angleA=90,angleB=-150](-.4,0)(0,.3)
\pccurve[doubleline=true, angleA=30,angleB=-120](0,.3)(.5,.7)
\pccurve[doubleline=true, angleA=60,angleB=-90](.5,.7)(.6,1.3)
\pccurve[doubleline=true, angleA=90,angleB=-60](.6,1.3)(.4,1.8)
\pccurve[doubleline=true, angleA=120,angleB=0](.4,1.8)(0,2)
\pccurve[doubleline=true, angleA=180,angleB=90](0,2)(-.8,1.45)
\pccurve[doubleline=true, angleA=-90,angleB=170](-.8,1.45)(-.6,1.1)
\pccurve[doubleline=true, angleA=-10,angleB=180](-.6,1.1)(0,1.05)
\pccurve[doubleline=true, angleA=0,angleB=-160](0,1.05)(.8,1.2)
\pccurve[doubleline=true, angleA=20,angleB=-90](.8,1.2)(1.1,1.5)
\pccurve[doubleline=true, angleA=90,angleB=0](1.1,1.5)(.9,1.8)
\pccurve[doubleline=true, angleA=180,angleB=50](.9,1.8)(.4,1.45)
\pccurve[doubleline=true, angleA=-130,angleB=70](.4,1.45)(.2,1)
\pccurve[doubleline=true, angleA=-110,angleB=135](.2,1)(.4,.6)
\pccurve[doubleline=true, angleA=-45,angleB=90](.4,.6)(1.2,0)
\pccurve[doubleline=true, angleA=-90,angleB=-45](1.2,0)(0,-.4)
\pccurve[doubleline=true, angleA=60,angleB=0](-.7,.8)(-.7,1.45)
\pccurve[doubleline=true, angleA=0,angleB=-160](0,1.05)(.8,1.2)
\pccurve[doubleline=true, angleA=60,angleB=-90](.5,.7)(.6,1.3)
\pscircle[linecolor=white, fillstyle=solid, fillcolor=lightgray](0,-.4){.3}
\psarc(0,-.4){.3}{-150}{-29} \psarc(0,-.4){.3}{-18}{30} \psarc(0,-.4){.3}{40}{140} \psarc(0,-.4){.3}{150}{198}
\rput(0,-1.5){{$(d)$}}
\rput(0,-.4){{$\mathcal{D}$}}
\end{pspicture}
$$
\caption{$(a)$ A $2$-band $(b)$ its banded surface with a blackboard framing, $(c)$ a bouquet of $2$-circles presentation of $(a)$
and $(d)$ its bouquet of circles surfaces with a blackboard framing.} \label{bandsurface}
\end{figure}
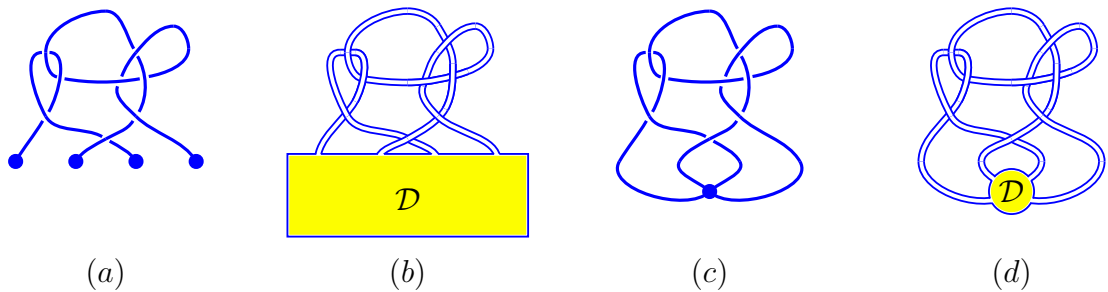

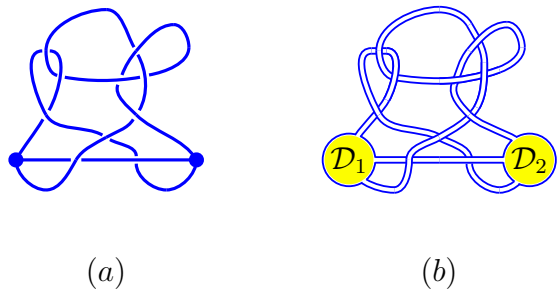
\begin{figure}
$$
\begin{pspicture}[shift=-1.3](-1.4,-1.7)(1.4,2.1)
\pscircle[fillcolor=lightgray, fillstyle=solid, linewidth=3pt](-1.2,0){.1}
\pscircle[fillcolor=lightgray, fillstyle=solid, linewidth=3pt](1.2,0){.1}
\pccurve[linewidth=1.2pt, angleA=60,angleB=-120](-1.2,0)(-.84,.54)
\pccurve[linewidth=1.2pt, angleA=60,angleB=0](-.76,.66)(-.75,1.45)
\pccurve[linewidth=1.2pt, angleA=180,angleB=120](-.85,1.45)(-.8,.6)
\pccurve[linewidth=1.2pt, angleA=-60,angleB=150](-.8,.6)(-.05,.33)
\pccurve[linewidth=1.2pt, angleA=-30,angleB=90](.05,.27)(.4,.05)
\pccurve[linewidth=1.2pt, angleA=-90,angleB=180](.4,-.05)(.8,-.4)
\pccurve[linewidth=1.2pt, angleA=0,angleB=-90](.8,-.4)(1.2,0)
\pccurve[linewidth=1.2pt, angleA=-120,angleB=0](-.4,0)(-.8,-.4)
\pccurve[linewidth=1.2pt, angleA=180,angleB=-90](-.8,-.4)(-1.2,0)
\pccurve[linewidth=1.2pt, angleA=60,angleB=-150](-.4,0)(0,.3)
\pccurve[linewidth=1.2pt, angleA=30,angleB=-120](0,.3)(.36,.54)
\pccurve[linewidth=1.2pt, angleA=60,angleB=-60](.44,.66)(.42,1.4)
\pccurve[linewidth=1.2pt, angleA=120,angleB=0](.37,1.5)(0,2)
\pccurve[linewidth=1.2pt, angleA=180,angleB=90](0,2)(-.8,1.45)
\pccurve[linewidth=1.2pt, angleA=-90,angleB=170](-.8,1.45)(-.65,1.125)
\pccurve[linewidth=1.2pt, angleA=-10,angleB=-170](-.55,1.1)(.45,1.1)
\pccurve[linewidth=1.2pt, angleA=10,angleB=-90](.57,1.13)(1.1,1.5)
\pccurve[linewidth=1.2pt, angleA=90,angleB=0](1.1,1.5)(.9,1.8)
\pccurve[linewidth=1.2pt, angleA=180,angleB=50](.9,1.8)(.4,1.45)
\pccurve[linewidth=1.2pt, angleA=-130,angleB=70](.4,1.45)(.21,1.13)
\pccurve[linewidth=1.2pt, angleA=-110,angleB=135](.155,1.0)(.4,.6)
\pccurve[linewidth=1.2pt, angleA=-45,angleB=120](.4,.6)(1.2,0)
\psline[linewidth=1.2pt](-1.2,0)(-.47,0)
\psline[linewidth=1.2pt](-.33,0)(1.2,0)
\rput(0,-1.5){{$(a)$}}
\end{pspicture} \quad\quad
\begin{pspicture}[shift=-1.3](-2.2,-1.7)(2.2,2.1)
\psline[doubleline=true](-1.2,0)(0,0)
\pccurve[doubleline=true, angleA=90,angleB=-120](-1.2,0)(-.7,.8)
\pccurve[doubleline=true, angleA=60,angleB=0](-.7,.8)(-.7,1.45)
\pccurve[doubleline=true, angleA=180,angleB=0](-.7,1.45)(-.9,1.45)
\pccurve[doubleline=true, angleA=180,angleB=120](-.9,1.45)(-.8,.6)
\pccurve[doubleline=true, angleA=-60,angleB=150](-.8,.6)(0,.3)
\pccurve[doubleline=true, angleA=-30,angleB=90](0,.3)(.4,0)
\pccurve[doubleline=true, angleA=-90,angleB=180](.4,-.05)(.8,-.4)
\pccurve[doubleline=true, angleA=0,angleB=-90](.8,-.4)(1.2,0)
\pccurve[doubleline=true, angleA=0,angleB=-90](-.6,-.4)(-.4,0)
\pccurve[doubleline=true, angleA=180,angleB=-90](-.6,-.4)(-1.2,0)
\pccurve[doubleline=true, angleA=90,angleB=-150](-.4,0)(0,.3)
\pccurve[doubleline=true, angleA=30,angleB=-120](0,.3)(.5,.7)
\pccurve[doubleline=true, angleA=60,angleB=-90](.5,.7)(.6,1.3)
\pccurve[doubleline=true, angleA=90,angleB=-60](.6,1.3)(.4,1.8)
\pccurve[doubleline=true, angleA=120,angleB=0](.4,1.8)(0,2)
\pccurve[doubleline=true, angleA=180,angleB=90](0,2)(-.8,1.45)
\pccurve[doubleline=true, angleA=-90,angleB=170](-.8,1.45)(-.6,1.1)
\pccurve[doubleline=true, angleA=-10,angleB=180](-.6,1.1)(0,1.05)
\pccurve[doubleline=true, angleA=0,angleB=-160](0,1.05)(.8,1.2)
\pccurve[doubleline=true, angleA=20,angleB=-90](.8,1.2)(1.1,1.5)
\pccurve[doubleline=true, angleA=90,angleB=0](1.1,1.5)(.9,1.8)
\pccurve[doubleline=true, angleA=180,angleB=50](.9,1.8)(.4,1.45)
\pccurve[doubleline=true, angleA=-130,angleB=70](.4,1.45)(.2,1)
\pccurve[doubleline=true, angleA=-110,angleB=135](.2,1)(.4,.6)
\pccurve[doubleline=true, angleA=-45,angleB=90](.4,.6)(1.2,0)
\pccurve[doubleline=true, angleA=60,angleB=0](-.7,.8)(-.7,1.45)
\pccurve[doubleline=true, angleA=0,angleB=-160](0,1.05)(.8,1.2)
\pccurve[doubleline=true, angleA=60,angleB=-90](.5,.7)(.6,1.3)
\psline[doubleline=true](1.2,0)(0,0)
\pscircle[linecolor=white, fillstyle=solid, fillcolor=lightgray](-1.2,0){.35}
\pscircle[linecolor=white, fillstyle=solid, fillcolor=lightgray](1.2,0){.35}
\psarc(-1.2,0){.35}{5}{65} \psarc(-1.2,0){.35}{74}{301} \psarc(-1.2,0){.35}{311}{355}
\psarc(1.2,0){.35}{245}{121} \psarc(1.2,0){.35}{131}{175} \psarc(1.2,0){.35}{185}{235}
\rput(0,-1.5){{$(b)$}}
\rput(-1.2,0){{$\mathcal{D}_1$}} \rput(1.2,0){{$\mathcal{D}_2$}}
\end{pspicture}
$$
\caption{$(a)$ a $3$-dipole presentation of $3$-bands and $(b)$ its $3$-dipole surfaces with a blackboard framing.} \label{dipolefig}
\end{figure}

For fixed types of graphs, considering the graph embedding as Seifert surfaces is not new.
Kauffman obtained a \emph{banded surface} $S$ from an $n$-band $B$ by attaching
disc $\D$ with a blackboard framing as depicted in Figure~\ref{bandsurface} $(b)$.
One may find that these are embeddings of a bouquet of circles, and we call these
banded surfaces \emph{bouquet of $n$-surfaces}. An example of a $3$-dipole and its flat $3$-dipole surface are given in Figure~\ref{dipolefig}.
In particular, if all bands in a bouquet of $n$-surface
are flat, then we call it a \emph{flat bouquet of $n$-surface}. Similarly, we define
\emph{$n$-dipole surfaces}, \emph{flat $n$-dipole surfaces}, \emph{$K_{n,m}$ surfaces}
and \emph{flat $K_{n,m}$ surfaces}. Let us remark that $1$-bands in these surfaces may be linked but not knotted.
The author, Kwon and Lee~\cite{KKL:string} showed the existence of a flat unknotted banded surface whose boundary is the given link. 
They also provided a few classification theorems and some upperbound on the minimal numbers of bands required to present the given link.
As we have mentioned and illustrated in Figure~\ref{bandsurface}, these flat unknotted banded surfaces with $n$ bands may be considered as
bouquet of $n$-surfaces.

\begin{thm} (\cite{KKL:string})
For a given link $L$, there exists a flat unknotted banded surface $F$ whose boundary is the link $L$.
\label{existencethm}
\end{thm}

Furthermore, if $1$-bands in an $n$-dipole surface are braided, it is called a \emph{braidzel surface} and
Nakamura showed that every link is a boundary of a braidzel surface~\cite{Nakamura:braidzel}. His result was
improved by Miura~\cite{Miura} that every link is a boundary of a flat braidzel surface.

Bands in a braidzel surface are allowed to be twisted for any integral number of times. We will only consider
an even number of twists. Similarly, a framed banded surface can be produced by replacing each arc by a band
(each arc in the middle of the band is called a \emph{core}), where a prefixed
framing on the band represents $m_i$ full twists; precisely,
$m_i$ is the linking number between a closed path $\alpha$, which is a path
sum of the \emph{core} of the $i$-th band of $S$, and any path joining both ends of the core in
$\D$ and its push up $\alpha^{+}$ towards to the positive normal direction
(as indicated ``$+$" in Figure~\ref{bandsurface} $(b)$).
The linking number discussed here
does not depend on the choice of path on $\D$.

\section{Dipole surfaces and $K_{2,n}$ surfaces} \label{dipole}

A motivation to consider the dipole surfaces initiated from an article by Baader~\cite{Baader:bipartite}
which introduced \emph{ribbon diagrams for strongly quasipositive links}
in order to show that every $(m, n)$ torus link is a boundary of a surface which
is obtained from the $0$ voltage assignment on all edges of the complete bipartite graphs $K_{m,n}$, where
the diagram of the complete bipartite graph $K_{m,n}$ is chosen to be in a very special form as explained the \emph{standard diagram} as depicted in Figure~\ref{k23torus} for $K_{2,3}$.

\begin{figure}
$$
\begin{pspicture}[shift=-1.4](-.2,-.2)(4.5,2.6)
\psline[doubleline=true](1,2)(0,0)(3,2)(2,0)(1,2)(4,0)(3,2)
\psline[doubleline=true](1.6,1.6)(3.4,.4)
\psline[doubleline=true](.6,.4)(1.8,1.2)
\psline[linecolor=white, linewidth=10pt](-.25,0)(1,2.5)
\psline(-.0405,0)(.99,2.064)
\psline[linecolor=white, linewidth=10pt](4.25,0)(3,2.5)
\psline(4.0405,0)(3.01,2.064)
\end{pspicture}
\quad \cong \quad
\begin{pspicture}[shift=-1.8](-2.3,-1.5)(2.2,1.8)
\pccurve[angleA=90,angleB=-160](-.6,-.5)(-.07,-.03)
\pccurve[angleA=-90,angleB=160](-.6,-.5)(0,-1)
\pccurve[angleA=-20,angleB=180](0,-1)(.8,-1.2)
\pccurve[angleA=0,angleB=-90](.8,-1.2)(1.8,0)
\pccurve[angleA=90,angleB=-10](1.8,0)(.65,.98)
\pccurve[angleA=-90,angleB=30](.6,1)(.07,.03)
\pccurve[angleA=0,angleB=90](0,1.8)(.6,1)
\pccurve[angleA=90,angleB=180](-.6,1.05)(0,1.8)
\pccurve[angleA=160,angleB=-90](0,0)(-.6,.95)
\pccurve[angleA=90,angleB=-20](.6,-.5)(0,0)
\pccurve[angleA=20,angleB=-90](.07,-.97)(.6,-.5)
\pccurve[angleA=-160,angleB=0](-.07,-1.03)(-.8,-1.2)
\pccurve[angleA=180,angleB=-90](-.8,-1.2)(-1.8,0)
\pccurve[angleA=90,angleB=190](-1.8,0)(-.6,1)
\pccurve[angleA=10,angleB=170](-.6,1)(.54,1.01)
\end{pspicture}
$$
\caption{A different diagram of the complete bipartite graph $K_{2,3}$ whose boundary is the figure-eight knot.} \label{k23figure8}
\end{figure}
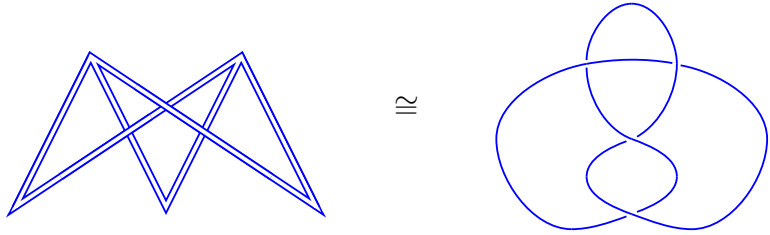

However, if we use a non-standard diagram of $K_{2,3}$, as depicted in Figure~\ref{k23figure8}, the boundary of this surface is the figure eight knot.
Since the figure eight knot is hyperbolic, it must not be a torus knot.
This phenomenon motivates us to find an answer for the question as to whether every link is a boundary of a complete bipartite graph $K_{n,m}$.
In~\cite{KKL:alternatingsign}, the following question was raised in order to settle this problem.

\begin{que} (\cite{KKL:alternatingsign}) \label{que1}
For a given link $L$, is there a graph diagram $D(\Gamma)$ of a complete bipartite graph $\Gamma$ such that the link $L$ is a boundary of $D(\Gamma)$ where
all voltage assignments on the edges of $\Gamma$ are $0$?
\end{que}

\begin{figure}
$$
\begin{pspicture}[shift=-1.4](-.3,-.3)(4.5,2.5)
\psline[doubleline=true](1,2)(0,0)(3,2)(2,0)(1,2)(4,0)(3,2)
\psline[doubleline=true](1.6,1.6)(3.4,.4)
\psline[doubleline=true](1.8,.4)(1.2,1.6)
\psline[linecolor=white, linewidth=10pt](-.25,0)(1,2.5)
\psline(-.0405,0)(.99,2.064)
\psline[linecolor=white, linewidth=10pt](4.25,0)(3,2.5)
\psline(4.0405,0)(3.01,2.064)
\end{pspicture}
\quad \cong \quad
\begin{pspicture}[shift=-.8](-2.3,-0.7)(2.2,1.3)
\psarc(.75,.25){1}{0}{90}
\psarc(.75,-.25){1}{270}{0}
\psarc(-.75,.25){1}{90}{180}
\psarc(-.75,-.25){1}{180}{270}
\psarc(.75,.25){.5}{0}{180}
\psarc(.75,-.25){.5}{180}{0}
\psarc(-.75,.25){.5}{0}{180}
\psarc(-.75,-.25){.5}{180}{0}
\psline(-.75,1.25)(.75,1.25)
\psline(-.75,-1.25)(.75,-1.25)
\pccurve[angleA=-90,angleB=90](-1.75,.25)(-1.25,-.25)
\pccurve[angleA=-90,angleB=90](-.25,.25)(.25,-.25)
\pccurve[angleA=-90,angleB=90](1.25,.25)(1.75,-.25)
\pccurve[angleA=-90,angleB=45](1.75,.25)(1.55,.05)
\pccurve[angleA=-135,angleB=90](1.45,-.05)(1.25,-.25)
\pccurve[angleA=-90,angleB=45](.25,.25)(.1,.05)
\pccurve[angleA=-135,angleB=90](-.1,-.05)(-.25,-.25)
\pccurve[angleA=-90,angleB=45](-1.25,.25)(-1.45,.05)
\pccurve[angleA=-135,angleB=90](-1.55,-.05)(-1.75,-.25)
\end{pspicture}
$$
\caption{A diagram of the complete bipartite graph $K_{2,3}$ whose boundary is the torus knot $T(2,3)$.} \label{k23torus}
\end{figure}
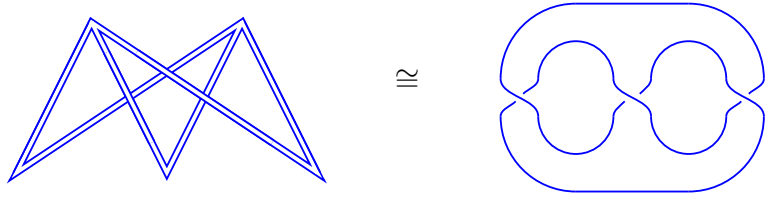

A weaker version of Question~\ref{que1} was proven by Kim and et al. in~\cite{JJK:K2n}.

\begin{thm} (\cite{JJK:K2n}) \label{main2}
For a given link $L$, there exists a graph diagram $D(K_{2,n})$ of a complete bipartite graph $K_{2,n}$ such that the link $L$ is a boundary of $D(K_{2,n})$ where
all voltage assignments on the edges of $K_{2,n}$ are either $0, 1$ or $-1$.
\end{thm}

Now, we will prove Theorem~\ref{dipolemainthm}.
First, we explain the moves which will be used in the proof of
the theorem. From a flat bouquet of $n$-circles surface, we can obtain a sequence of length $2n$ which represents
the connection of bands as follows. The disc $\mathcal{D}$ can be presented as unit disc,
the intervals which are the intersection of bands and the disc $\mathcal{D}$ can be labeled by
$\{ 1, 2$, $\ldots$, $n$, $\overline{1}$, $\overline{2}$, $\ldots$, $\overline{n} \}$ clockwisely
from a fixed point on the boundary of the disc $\mathcal{D}$ which is not in the intervals.
By moving the fixed point, we will get a different labels on the intervals, this is called
a \emph{relabeling}.
Then each bands connecting two interval can be presented by a subset of labels of the cardinality $2$.
The set of all such subsets is denoted by $L(\mathcal{F})$ from a flat banded surface with $n$ bands.
For example, $L(\mathcal{F})$ for a flat banded surface with $5$ bands in Figure~\ref{412complete} $(b)$
is $L(\mathcal{F}) = \{\{1, \overline{1}\}$, $\{2, \overline{4}\}$, $\{3, \overline{5}\}$, $\{4, \overline{2}\}$,
$\{5, \overline{3}\}\}$. Although all $i$ is paired with $\overline{j}$ in this example,
this is not true in general. Indeed, showing such $L(\mathcal{F})$ can
be chosen by a suitable relabeling and a finite sequence of band slides is the key idea of Theorem~\ref{dipolemainthm}.

\begin{thm} \label{dipolemainthm}
Every link is a boundary of an $n$-dipole flat surface.
\begin{proof}
For a given link $L$, there exists a flat bouquet of $n$-circles surface $B$ such that the boundary of $B$ is $L$ by Theorem~\ref{existencethm}.
From a flat bouquet of $n$-circles surface $B$ such that the boundary of $B$ is $L$,
as described the above, by fixing a point $P$ on the boundary of $\D$, we obtain a sequence $L(B)$.
First we prove the following claim.
\vskip .3cm
{\bf Claim} : any flat bouquet of $n$-circles surface $B$ can be transformed to a flat bouquet of $n$-circles surface $\mathcal{F}$ by
a finite sequence of band slides such that $L(\mathcal{F})$ consists only of the shape
$\{i, \overline{j}\}$. 
\vskip .3cm
To prove the claim we induct on $n$.

\begin{figure}
$$
\begin{pspicture}[shift=-1.2](-1.4,-2.1)(11.2,3)
\psarc[doubleline=true](2.1,0){.7}{-5}{185}
\psarc[doubleline=true](5.7,0){2.5}{-5}{185}
\psarc[doubleline=true](6.1,0){.7}{-5}{185}
\psline[doubleline=true](0,0)(0,.4)
\psline[doubleline=true](4.6,0)(4.6,.4)
\psline[doubleline=true](9.6,0)(9.6,.4)
\psframe[linecolor=lightgray,fillstyle=solid,fillcolor=lightgray](-1.2,-1.5)(10.8,0)
\psline(-1.2,0)(-.03,0) \psline(.03,0)(.4,0) \psline(1,0)(1.37,0)
\psline(1.43,0)(1.8,0) \psline(2.4,0)(2.77,0) \psline(2.83,0)(3.17,0)
\psline(3.23,0)(3.6,0) \psline(4.2,0)(4.57,0) \psline(4.63,0)(5.37,0)
\psline(5.43,0)(5.8,0) \psline(6.4,0)(6.77,0) \psline(6.83,0)(7.2,0)
\psline(7.8,0)(8.17,0) \psline(8.23,0)(8.6,0) \psline(9.2,0)(9.57,0)
\psline(9.63,0)(10.8,0)(10.8,-1.5)(-1.2,-1.5)(-1.2,0)
\rput(0,-.4){{$1$}} \rput(1.4,-.4){{$i$}} \rput(2.7,-.4){{$j$}}
\rput(3.3,-.4){{$j+1$}} \rput(4.6,-.4){{$n$}} \rput(5.4,-.4){{$\overline{1}$}}
\rput(6.8,-.4){{$\overline{k}$}} \rput(8.2,-.4){{$\overline{l}$}} \rput(9.6,-.4){{$\overline{n}$}}
\rput(.7,0){{$\ldots$}} \rput(2.1,0){{$\ldots$}} \rput(3.9,0){{$\ldots$}}
\rput(6.1,0){{$\ldots$}} \rput(7.5,0){{$\ldots$}} \rput(8.9,0){{$\ldots$}}
\psline[linecolor=darkred,linewidth=1.5pt, linestyle=dashed](5,-.5)(5,2.8)
\rput(5,-.9){{$\mathcal{D}$}}
\rput(5,-1.9){{$(a)$}}
\end{pspicture}
$$
$$
\begin{pspicture}[shift=-1.2](-1.4,-2.1)(11.2,3.5)
\psarc[doubleline=true](2.1,0){.7}{90}{185}
\pccurve[doubleline=true, angleA=0,angleB=-120](2.1,.7)(2.817,1.55)
\psarc[doubleline=true](5.5,0){3.1}{-5}{150}
\psarc[doubleline=true](5.5,0){2.7}{-5}{185}
\psarc[doubleline=true](5.7,0){1.1}{-5}{185}
\psline[doubleline=true](0,0)(0,.4)
\psline[doubleline=true](4.2,0)(4.2,.4)
\psline[doubleline=true](5.4,0)(5.4,.4)
\psline[doubleline=true](10,0)(10,.4)
\psframe[linecolor=lightgray,fillstyle=solid,fillcolor=lightgray](-1.2,-1.5)(10.8,0)
\psline(-1.2,0)(-.03,0) \psline(.03,0)(.4,0) \psline(1,0)(1.37,0)
\psline(1.43,0)(1.8,0) \psline(2.4,0)(2.77,0) \psline(2.83,0)(3.2,0)
\psline(3.8,0)(4.17,0) \psline(4.23,0)(4.57,0) \psline(4.63,0)(5.37,0)
\psline(5.43,0)(5.8,0) \psline(6.4,0)(6.77,0) \psline(6.83,0)(7.2,0)
\psline(7.8,0)(8.17,0) \psline(8.23,0)(8.57,0) \psline(8.63,0)(9,0)
\psline(9.6,0)(9.97,0) \psline(10.03,0)(10.8,0)(10.8,-1.5)(-1.2,-1.5)(-1.2,0)
\rput(0,-.4){{$1$}} \rput(1.4,-.4){{$i$}} \rput(2.8,-.4){{$j$}}
\rput(4.6,-.4){{$n$}} \rput(5.4,-.4){{$\overline{1}$}}
\rput(6.8,-.4){{$\overline{k-1}$}} \rput(8.0,-.4){{$\overline{l-1}$}}
\rput(8.6,-.4){{$\overline{l}$}} \rput(10,-.4){{$\overline{n}$}}
\rput(.7,0){{$\ldots$}} \rput(2.1,0){{$\ldots$}} \rput(3.5,0){{$\ldots$}}
\rput(6.1,0){{$\ldots$}} \rput(7.5,0){{$\ldots$}} \rput(9.8,0){{$\ldots$}}
\psline[linecolor=darkred,linewidth=1.5pt, linestyle=dashed](5,-.5)(5,3.3)
\rput(5,-.9){{$\mathcal{D}$}}
\rput(5,-1.9){{$(b)$}}
\end{pspicture}
$$
\caption{A slide of the band connects $i$ and $j$ along the band $j+1$ and $\overline{l}$
if $1 \le i < j < n < \overline{1}< \overline{k},  \overline{l} \le \overline{n}$.} \label{lempf1}
\end{figure}
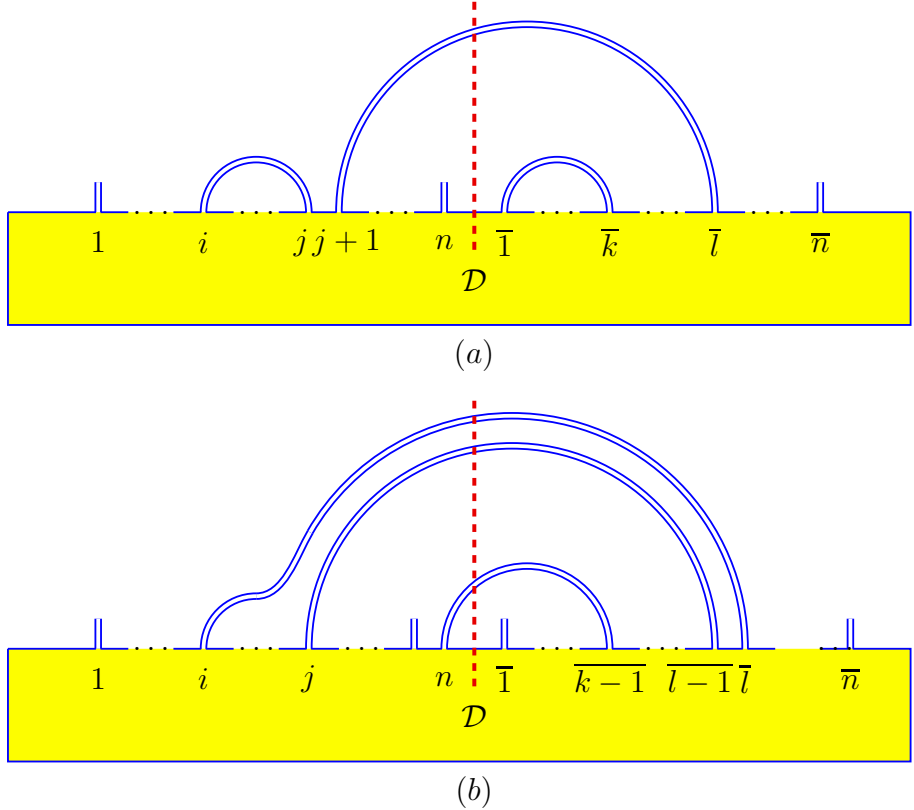

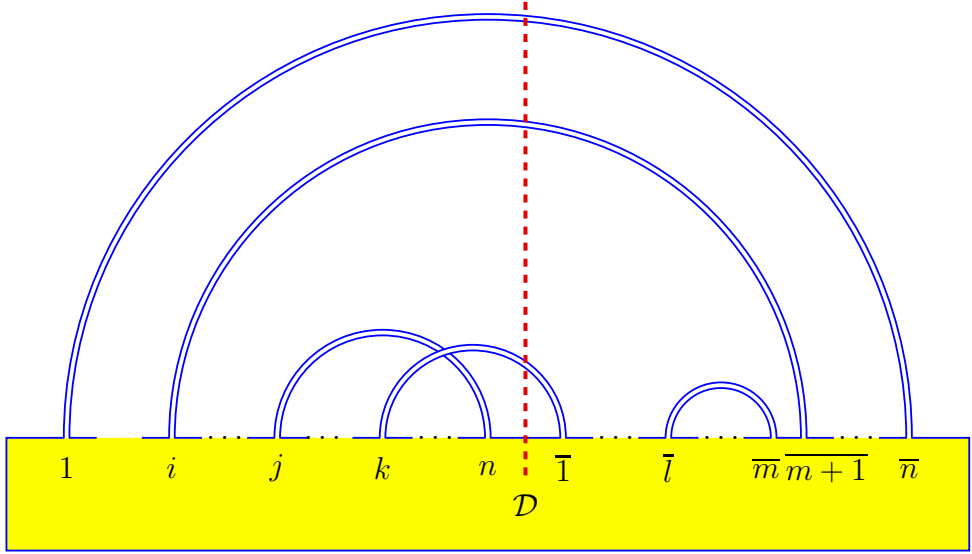
\begin{figure}
$$
\begin{pspicture}[shift=-1.2](-2.8,-2.1)(10.4,6)
\psarc[doubleline=true](2.4,0){1.4}{-5}{185}
\psarc[doubleline=true](3.6,0){1.2}{-5}{185}
\psarc[doubleline=true](3.8,0){4.2}{-5}{185}
\psarc[doubleline=true](3.8,0){5.6}{-5}{185}
\psarc[doubleline=true](6.9,0){0.7}{-5}{185}
\psframe[linecolor=lightgray,fillstyle=solid,fillcolor=lightgray](-2.6,-1.5)(10.2,0)
\psline(-1.77,0)(-1.4,0) \psline(-.8,0)(-0.43,0) \psline(-0.37,0)(0,0) \psline(0.6,0)(.97,0)
\psline(1.03,0)(1.4,0) \psline(2,0)(2.37,0) \psline(2.43,0)(2.8,0)
\psline(3.4,0)(3.77,0) \psline(3.83,0)(4.77,0) \psline(4.83,0)(5.2,0)
\psline(5.8,0)(6.17,0) \psline(6.23,0)(6.6,0) \psline(7.2,0)(7.57,0) \psline(7.63,0)(7.97,0)
\psline(8.03,0)(8.4,0) \psline(9,0)(9.37,0) \psline(9.43,0)(10.2,0)(10.2,-1.5)(-2.6,-1.5)(-2.6,0)(-1.83,0)
\rput(-1.8,-.4){{$1$}} \rput(-.4,-.4){{$i$}} \rput(1,-.4){{$j$}} \rput(2.4,-.4){{$k$}}
\rput(3.8,-.4){{$n$}} \rput(4.8,-.4){{$\overline{1}$}}
\rput(6.2,-.4){{$\overline{l}$}} \rput(7.5,-.4){{$\overline{m}$}}
\rput(8.3,-.4){{$\overline{m+1}$}} \rput(9.4,-.4){{$\overline{n}$}}
\rput(.3,0){{$\ldots$}} \rput(1.6,0){{$\ldots$}} \rput(3.1,0){{$\ldots$}}
\rput(5.5,0){{$\ldots$}} \rput(6.9,0){{$\ldots$}} \rput(8.7,0){{$\ldots$}}
\psline[linecolor=darkred,linewidth=1.5pt, linestyle=dashed](4.3,-.5)(4.3,5.8)
\rput(4.3,-.9){{$\mathcal{D}$}}
\end{pspicture}
$$
\caption{A slide of the band connects $i$ and $j$ along the band $j+1$ and $\overline{l}$
if $1 \le i < $$j < n < $$\overline{1} < $$\overline{k} < $$ \overline{l} $$\le \overline{n}$.} \label{lempf2}
\end{figure}

If $n=1$, $L(\F)$ is clearly the desired form. For $n=2$, there are three possible cases for $L(\mathcal{F})$ :
$\{\{1, 2\}$, $\{ \overline{1}, \overline{2}\}\}$, $\{\{1, \overline{1}\}$, $\{ 2, \overline{2}\}\}$
and $\{\{1, \overline{2}\}$, $\{ 2, \overline{1}\}\}$, but the last two of them are already in the desired form.
For $\{\{1, 2\}$, $\{ \overline{1}, \overline{2}\}\}$, we slide the first bands presented by $\{1, 2\}$
along the band presented by $\{ \overline{1}, \overline{2}\}\}$ to have $\{\{1, \overline{2}\}$, $\{ 2, \overline{1}\}\}$.

Now we assume $n \ge 3$. Suppose there exist a flat bouquet of $n$-circles surface $B$ which
can not be transformed to a flat bouquet of $n$-circles surface $\mathcal{F}$ by
a finite sequence of band slides such that $L(\mathcal{F})$ consists only of the shape
$\{i, \overline{j}\}$ for a suitable fixed point which determines the relabeling.
Then, for any fixed point $R$ on the boundary of the disc $\mathcal{D}$ and $\mathcal{F}$ which is obtained
from $\mathcal{F}$ by a finite sequence of band slides, the set $\Omega(\mathcal{F}, R) =$
$\{\{ i, j\} \in L(\mathcal{F}) | 1 \le i< j \le n\}$ is nonempty.

We pick one $\mathcal{F}_0$ which has the minimal cardinality of the set $\Omega(\mathcal{F}_0, P)$ for a point $P$ among all counterexamples.
We choose a pair $\{i,j\}$ for which $j$ is the largest among such pairs in $\Omega(\mathcal{F}_0, P)$.
If $j \not= n$, then by the maximality of $j$, $j+1$ (possibly $j+1$ can be $n$) must be connected to $\overline{l}$.
We divide cases $1-i)$ $\overline{1}$ is connected to $\overline{k}$ and $1-ii)$ $\overline{1}$ is connected to $k$.

For the first case $1-i)$ as shown in Figure~\ref{lempf1}, we slide the band presented by $\{ i, j\}$ along the band presented by $\{ j+1, \overline{l}\}$.
The resulting flat bouquet of $n$-circles surface $\mathcal{F}_1$ has a new $L(\mathcal{F},P)$
which has the property that $\Omega(\mathcal{F}_1, P) < \Omega(\mathcal{F}_0, P)$ where we used the same fixed point $P$ for labeling.
But it contradicts the minimum hypothesis of $\Omega(\mathcal{F}_0, P)$.

For the second case $1-ii)$, even if we slide the band presented by $\{ i, j\}$ along the band presented by $\{ j+1, \overline{l}\}$,
the resulting flat bouquet of $n$-circles surface $\mathcal{F}_1$ has the property that
$\Omega(\mathcal{F}_1, P) = \Omega(\mathcal{F}_0, P)$ but
$\{k,n\} \in \Omega(\mathcal{F}_1, P)$. Then this case can be handled as follows.

Now without the loss of generality, we may assume that there exists a pair $\{i,n\}$ in $\Omega(\mathcal{F}_0, P)$.
Similarly, we divide cases $2-i)$ $\overline{1}$ is connected to $\overline{k}$ and $2-ii)$ $\overline{1}$ is connected to $k$.

For the first case $2-i)$, if $\overline{1}$ is connected to $\overline{k}$, then we slide the band presented by $\{ i, n\}$ along the band presented by $\{ \overline{1}, \overline{k}\}$
and it leads us the similar contradiction for the case $1-i)$.

For the second case $2-ii)$ if $\overline{1}$ is connected to $k$, then we move $P$ to a point $Q$ between $\overline{n-1}$ and $\overline{n}$.
Then $\Omega(\mathcal{F}_0, P) \ge \Omega(\mathcal{F}_0, Q)$ but
by the minimum hypothesis of $\Omega(\mathcal{F}_0, P)$, we must have $\Omega(\mathcal{F}_0, P) = \Omega(\mathcal{F}_0, Q)$.
To obtain such an equality, one can see that $\{1, \overline{n}\} \in \Omega(\mathcal{F}_0, P)$. Now, we consider
the set $\overline{\Omega}(\mathcal{F}_0, P) =$ $\{\{ \overline{i}, \overline{j}\} \in L(\mathcal{F}) | 1 \le i< j \le n \}$.
One can see that the cardinality of two sets $\overline{\Omega}(\mathcal{F}, P)$, $\Omega(\mathcal{F}, P)$ must be the same. 
Therefore, there exists
$\{\overline{l}, \overline{m}\} \in \overline{\Omega}(\mathcal{F}, P)$ and we further choose one $\{\overline{l}, \overline{m}\}$
for which $m$ is the largest among such pairs in $\overline{\Omega}(\mathcal{F}_0, P)$ as illustrated in Figure~\ref{lempf2}.
If we slide the band presented by $\{\overline{l}, \overline{m}\}$ along the band presented by $\{ i, \overline{m+1}\}$.
The resulting flat bouquet of $n$-circles surface $\mathcal{F}_1$ has a new $L(\mathcal{F}_1)$ with the same fixed point $P$
which has the property that $\Omega(\mathcal{F}_1, P) < \Omega(\mathcal{F}_0, P)$.
It contradicts the minimum hypothesis of  $\Omega(\mathcal{F}_0, P)$ and completes the proof of the claim.

\begin{figure}
$$
\begin{pspicture}[shift=-2.1](-2.3,-2.7)(2.3,2.2)
\pscircle[linecolor=lightgray, fillcolor=lightgray, fillstyle=solid](0,0){2}
\pccurve[doubleline=true, angleA=-40,angleB=180](2;140)(-.6,1)
\pccurve[doubleline=true, angleA=-20,angleB=180](2;160)(-.6,.6)
\pccurve[doubleline=true, angleA=0,angleB=180](2;180)(-.6,.2)
\pccurve[doubleline=true, angleA=40,angleB=180](2;220)(-.6,-1)
\pccurve[doubleline=true, angleA=-140,angleB=0](2;40)(.6,1)
\pccurve[doubleline=true, angleA=-160,angleB=0](2;20)(.6,.6)
\pccurve[doubleline=true, angleA=180,angleB=0](2;0)(.6,.2)
\pccurve[doubleline=true, angleA=140,angleB=0](2;-40)(.6,-1)
\psarc(0,0){2}{1}{19} \psarc(0,0){2}{21}{39}
\psarc(0,0){2}{41}{139} \psarc(0,0){2}{141}{159} \psarc(0,0){2}{161}{179}
\psarc(0,0){2}{181}{190} \psarc(0,0){2}{210}{219} \psarc(0,0){2}{221}{319}
\psarc(0,0){2}{321}{330} \psarc(0,0){2}{350}{359}
\rput(0,1){{$\ldots$}} \rput(0,.6){{$\ldots$}}
\rput(0,.2){{$\ldots$}} \rput(0,-1){{$\ldots$}}
\rput(-1.8,-.6){{$\ldots$}} \rput(1.8,-.6){{$\ldots$}}
\rput(0,-.5){{$\mathcal{D}$}}
\rput(0,-2.5){{$(a)$}}
\end{pspicture} \quad \Rightarrow \quad
\begin{pspicture}[shift=-2.1](-2.3,-2.7)(2.3,2.2)
\pscircle[linecolor=lightgray, fillcolor=lightgray, fillstyle=solid](0,0){2}
\psframe[linecolor=white, fillcolor=white, fillstyle=solid](-1,-2.1)(1,2.1)
\pscircle[linecolor=lightgray, fillcolor=lightgray, fillstyle=solid](1.8;120){.2}
\pscircle[linecolor=lightgray, fillcolor=lightgray, fillstyle=solid](1.8;60){.2}
\pscircle[linecolor=lightgray, fillcolor=lightgray, fillstyle=solid](1.8;240){.2}
\pscircle[linecolor=lightgray, fillcolor=lightgray, fillstyle=solid](1.8;-60){.2}
\psframe[linecolor=lightgray, fillcolor=lightgray, fillstyle=solid](-1,-1.5588)(-.7,1.5588)
\psframe[linecolor=lightgray, fillcolor=lightgray, fillstyle=solid](1,-1.5588)(.7,1.5588)
\psframe[linecolor=lightgray, fillcolor=lightgray, fillstyle=solid](-.7,-.1)(.7,.1)
\psline(-.7,1.5588)(-.7,.1)(.7,.1)(.7,1.5588) \psline(.7,-1.5588)(.7,-.1)(-.7,-.1)(-.7,-1.5588)
\pccurve[doubleline=true, angleA=-40,angleB=180](2;140)(-.6,1)
\pccurve[doubleline=true, angleA=-20,angleB=180](2;160)(-.6,.6)
\pccurve[doubleline=true, angleA=0,angleB=180](2;180)(-.6,.2)
\pccurve[doubleline=true, angleA=40,angleB=180](2;220)(-.6,-1)
\pccurve[doubleline=true, angleA=-140,angleB=0](2;40)(.6,1)
\pccurve[doubleline=true, angleA=-160,angleB=0](2;20)(.6,.6)
\pccurve[doubleline=true, angleA=180,angleB=0](2;0)(.6,.2)
\pccurve[doubleline=true, angleA=140,angleB=0](2;320)(.6,-1)
\psarc(0,0){2}{1}{19} \psarc(0,0){2}{21}{39}
\psarc(0,0){2}{41}{60} \psarc(0,0){2}{120}{139}
\psarc(0,0){2}{141}{159} \psarc(0,0){2}{161}{179}
\psarc(0,0){2}{181}{190} \psarc(0,0){2}{210}{219}
\psarc(0,0){2}{221}{240} \psarc(0,0){2}{300}{319}
\psarc(0,0){2}{321}{330} \psarc(0,0){2}{350}{359}
\psarc(1.8;60){.2}{60}{180}
\psarc(1.8;120){.2}{0}{120}
\psarc(1.8;240){.2}{240}{0}
\psarc(1.8;300){.2}{180}{300}
\rput(0,1){{$\ldots$}} \rput(0,.6){{$\ldots$}}
\rput(0,.2){{$\ldots$}} \rput(0,-1){{$\ldots$}}
\rput(-1.8,-.6){{$\ldots$}} \rput(1.8,-.6){{$\ldots$}}
\rput(-1.2,-.5){{$\mathcal{D}_1$}}
\rput(1.2,-.5){{$\mathcal{D}_2$}}
\rput(0,-2.5){{$(b)$}}
\end{pspicture}
$$
\caption{$(a)$ The resulting flat bouquet of $n$-circles surface $\overline{B}$ by moves in the claim and $(b)$ the final flat dipole surface.} \label{dipolefig1}
\end{figure}
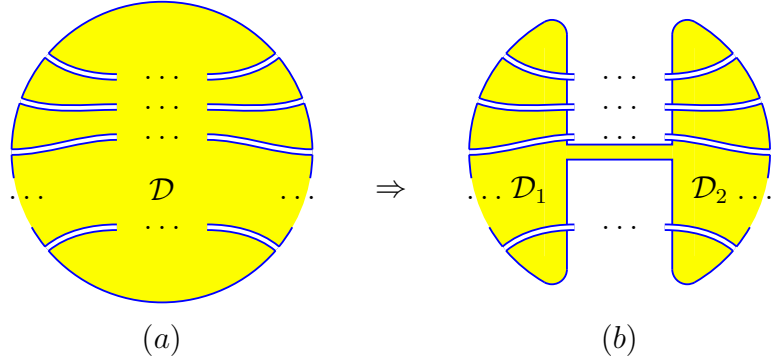

By the claim, the resulting flat bouquet of $n$-circles surface $\overline{B}$ has a shape as illustrated in Figure~\ref{dipolefig1} $(a)$.
The final modification is to make the disc $\D$ into two discs connected by a single flat band as depicted in Figure~\ref{dipolefig1} $(b)$.
The moves used in this process do not change the link type of the boundaries of the surfaces. This completes the proof of the theorem.
\end{proof}
\end{thm}

Let us remark that every braidzel surface is a dipole surface; however, the converse is not true in general.
However, the dipole surfaces dealt in the article, which are obtained from a flat plumbing basket surface are braidzel surfaces.
Moreover, they are completely determined by two permutations on $n$ letters, describing how bands are connected and the order that explains the layer of bands.

The following corollary positively answer the Question~\ref{que1}.

\begin{cor} \label{dipolemaincor}
For a given link $L$, there exists a graph diagram $D(K_{2,n})$ of a complete bipartite graph $K_{2,n}$ such that link $L$ is a boundary of $D(K_{2,n})$ where
all voltage assignments on the edges of $K_{2,n}$ are $0$.
\begin{proof}
For a given link $L$, there exists an $n$-dipole flat surface $B$ whose boundary is $L$ by Theorem~\ref{dipolemainthm}.
By subdividing each band in $B$ into two flat bands, the resulting surface can be considered a complete bipartite graph surface, where the bipartition
of the vertex set is (the set of vertices presented by two vertices in $n$-dipole, the set of vertices produced by the subdivision).
\end{proof}
\end{cor}

Let us remark that by flipping two discs corresponding to vertices in one of the bipartition sets whose cardinality is $2$, Corollary~\ref{dipolemaincor} can be restated that
for a given link $L$, there exists a graph diagram $D(K_{2,n})$ of a complete bipartite graph $K_{2,n}$ such that the link $L$ is a boundary of $D(K_{2,n})$ where
all voltage assignments on the edges of $K_{2,n}$ are $\pm1$.

\section{Algorithm and examples} \label{examples}

In this section, for a given link we provide how to find a graph diagram $D(K_{2,n})$ of a complete bipartite graph $K_{2,n}$ such that link $L$ is a boundary of $D(K_{2,n})$.
We also demonstrate our algorithm by exhibiting two examples, links $5_2$ and $4_1^2$.

\begin{figure}
$$
\begin{pspicture}[shift=-1.2](-.2,-1.9)(6,2)
\psline(0,0)(0,1)
\pccurve[angleA=90,angleB=180](0,1)(1,1.8)
\psline(1,1.8)(1.3,1.8)
\pccurve[angleA=0,angleB=180](1.3,1.8)(2.1,1)
\pccurve[angleA=0,angleB=-135](2.1,1)(2.45,1.3)
\pccurve[angleA=45,angleB=180](2.55,1.5)(2.9,1.8)
\pccurve[angleA=0,angleB=180](2.9,1.8)(3.7,1)
\pccurve[angleA=0,angleB=-135](3.7,1)(4,1.3)
\pccurve[angleA=45,angleB=180](4.2,1.5)(4.5,1.8)
\psline(4.5,1.8)(4.8,1.8)
\pccurve[angleA=0,angleB=90](4.8,1.8)(5.8,1)
\psline(5.8,1)(5.8,0)
\pccurve[angleA=-90,angleB=0](5.8,0)(4.8,-.8)
\psline(4.8,-.8)(1,-.8)
\pccurve[angleA=180,angleB=-90](1,-.8)(0,0)
\pccurve[angleA=180,angleB=180](1,0)(1,1)
\psline(1,1)(1.3,1)
\pccurve[angleA=0,angleB=-135](1.3,1)(1.6,1.3)
\pccurve[angleA=45,angleB=180](1.8,1.5)(2.1,1.8)
\pccurve[angleA=0,angleB=180](2.1,1.8)(2.9,1)
\pccurve[angleA=0,angleB=-135](2.9,1)(3.2,1.3)
\pccurve[angleA=45,angleB=180](3.4,1.5)(3.7,1.8)
\pccurve[angleA=0,angleB=180](3.7,1.8)(4.5,1)
\psline(4.5,1)(4.8,1)
\pccurve[angleA=0,angleB=0](4.8,1)(4.8,0)
\psline(4.8,0)(1,0)
\rput(2.9,-1.7){{$(a)$}}
\end{pspicture} \quad
\begin{pspicture}[shift=-1.2](-.7,-1.9)(5.2,2)
\psarc[doubleline=true](2.75,0){.75}{-5}{185}
\psarc[doubleline=true](2.25,0){.75}{-5}{185}
\psarc[doubleline=true](2.75,0){1.75}{-5}{185}
\psarc[doubleline=true](2.25,0){1.75}{-5}{185}
\psarc[doubleline=true](1.25,0){1.25}{-5}{185}
\psframe[linecolor=lightgray,fillstyle=solid,fillcolor=lightgray](-.5,-1)(5,0)
\psline(-.03,0)(-.5,0)(-.5,-1)(5,-1)(5,0)(4.53,0)
\psline(.03,0)(.47,0) \psline(.53,0)(.97,0) \psline(1.03,0)(1.47,0)
\psline(1.53,0)(1.97,0) \psline(2.03,0)(2.47,0) \psline(2.53,0)(2.97,0)
\psline(3.03,0)(3.47,0) \psline(3.53,0)(3.97,0) \psline(4.03,0)(4.47,0)
\rput(2.25,-.5){{$\mathcal{D}$}}
\rput(1.75,-1.7){{$(b)$}}
\end{pspicture}
$$

$$
\begin{pspicture}[shift=-.8](-2.8,-3.1)(2.8,2.2)
\pscircle[linecolor=lightgray, fillcolor=lightgray, fillstyle=solid](0,0){2}
\psframe[linecolor=white, fillcolor=white, fillstyle=solid](-1,-2.1)(1,2.1)
\pscircle[linecolor=lightgray, fillcolor=lightgray, fillstyle=solid](1.8;120){.2}
\pscircle[linecolor=lightgray, fillcolor=lightgray, fillstyle=solid](1.8;60){.2}
\pscircle[linecolor=lightgray, fillcolor=lightgray, fillstyle=solid](1.8;240){.2}
\pscircle[linecolor=lightgray, fillcolor=lightgray, fillstyle=solid](1.8;-60){.2}
\psframe[linecolor=lightgray, fillcolor=lightgray, fillstyle=solid](-1,-1.5588)(-.7,1.5588)
\psframe[linecolor=lightgray, fillcolor=lightgray, fillstyle=solid](1,-1.5588)(.7,1.5588)
\psframe[linecolor=lightgray, fillcolor=lightgray, fillstyle=solid](-.7,-.1)(.7,.1)
\psline(-.7,1.5588)(-.7,.1)(.7,.1)(.7,1.5588) \psline(.7,-1.5588)(.7,-.1)(-.7,-.1)(-.7,-1.5588)
\pccurve[doubleline=true, angleA=0,angleB=220](2;180)(2;40)
\pccurve[doubleline=true, angleA=40,angleB=180](2;220)(2;0)
\pccurve[doubleline=true, angleA=-20,angleB=200](2;160)(2;20)
\pccurve[doubleline=true, angleA=20,angleB=160](2;200)(2;-20)
\pccurve[doubleline=true, angleA=-40,angleB=140](2;140)(2;-40)
\psarc(0,0){2}{1}{19} \psarc(0,0){2}{21}{39}
\psarc(0,0){2}{41}{60} \psarc(0,0){2}{120}{139}
\psarc(0,0){2}{141}{159} \psarc(0,0){2}{161}{179}
\psarc(0,0){2}{181}{199} \psarc(0,0){2}{201}{219}
\psarc(0,0){2}{221}{240} \psarc(0,0){2}{300}{319}
\psarc(0,0){2}{321}{339} \psarc(0,0){2}{341}{359}
\psarc(1.8;60){.2}{60}{180}
\psarc(1.8;120){.2}{0}{120}
\psarc(1.8;240){.2}{240}{0}
\psarc(1.8;300){.2}{180}{300}
\rput(-1,1.3){{$\mathcal{D}_1$}}
\rput(1,1.3){{$\mathcal{D}_2$}}
\rput(2.3;140){{$1$}}
\rput(2.3;160){{$2$}}
\rput(2.3;180){{$3$}}
\rput(2.3;200){{$4$}}
\rput(2.3;220){{$5$}}
\rput(0,-3){{$(c)$}}
\end{pspicture} \quad
\begin{pspicture}[shift=-1.2](-.8,-1.2)(5.6,4.6)
\pccurve[doubleline=true, angleA=120,angleB=-140](5,0)(.5,3.7)
\pccurve[doubleline=true, angleA=60,angleB=-80](5,0)(4.1,3.5)
\pccurve[doubleline=true, angleA=120,angleB=-120](4,0)(.7,3.6)
\pccurve[doubleline=true, angleA=60,angleB=-60](4,0)(4.3,3.6)
\pccurve[doubleline=true, angleA=120,angleB=-100](3,0)(.9,3.5)
\pccurve[doubleline=true, angleA=60,angleB=-120](3,0)(3.7,3.6)
\pccurve[doubleline=true, angleA=120,angleB=-80](2,0)(1.1,3.5)
\pccurve[doubleline=true, angleA=60,angleB=-100](2,0)(3.9,3.5)
\pccurve[doubleline=true, angleA=120,angleB=-60](1,0)(1.3,3.6)
\pccurve[doubleline=true, angleA=60,angleB=-20](1,0)(4.5,3.7)
\pccurve[doubleline=true, angleA=120,angleB=-20](0,0)(1.5,3.7)
\pccurve[doubleline=true, angleA=60,angleB=-140](0,0)(3.5,3.7)
\pscircle[linecolor=lightgray, fillcolor=lightgray, fillstyle=solid](0,0){.3}
\pscircle[linecolor=lightgray, fillcolor=lightgray, fillstyle=solid](1,0){.3}
\pscircle[linecolor=lightgray, fillcolor=lightgray, fillstyle=solid](2,0){.3}
\pscircle[linecolor=lightgray, fillcolor=lightgray, fillstyle=solid](3,0){.3}
\pscircle[linecolor=lightgray, fillcolor=lightgray, fillstyle=solid](4,0){.3}
\pscircle[linecolor=lightgray, fillcolor=lightgray, fillstyle=solid](5,0){.3}
\pscircle[linecolor=lightgray, fillcolor=lightgray, fillstyle=solid](1,4){.6}
\pscircle[linecolor=lightgray, fillcolor=lightgray, fillstyle=solid](4,4){.6}
\psarc(0,0){.3}{112}{53} \psarc(0,0){.3}{63}{101}
\psarc(1,0){.3}{116.5}{50} \psarc(1,0){.3}{60}{106}
\psarc(2,0){.3}{123}{53} \psarc(2,0){.3}{63.5}{113}
\psarc(3,0){.3}{128}{59} \psarc(3,0){.3}{69}{117}
\psarc(4,0){.3}{130}{56} \psarc(4,0){.3}{67}{120}
\psarc(5,0){.3}{130}{61} \psarc(5,0){.3}{73}{120}
\psarc(1,4){.6}{-28}{209} \psarc(1,4){.6}{214}{232}
\psarc(1,4){.6}{237.5}{256.5} \psarc(1,4){.6}{262}{278.5}
\psarc(1,4){.6}{284}{302.5} \psarc(1,4){.6}{308}{327}
\psarc(4,4){.6}{-28}{209} \psarc(4,4){.6}{214}{232}
\psarc(4,4){.6}{237.5}{256.5} \psarc(4,4){.6}{262}{278.5}
\psarc(4,4){.6}{284}{302.5} \psarc(4,4){.6}{308}{327}
\rput(2.5,-.8){{$(d)$}}
\rput(1,4){{$\mathcal{D}_1$}}
\rput(4,4){{$\mathcal{D}_2$}}
\rput(0,0){{$0$}}
\rput(1,0){{$1$}}
\rput(2,0){{$2$}}
\rput(3,0){{$3$}}
\rput(4,0){{$4$}}
\rput(5,0){{$5$}}
\end{pspicture}
$$
\caption{$(a)$ The link $4_1^2$, $(b)$ a flat plumbing basket surface of the link $4_1^2$ with $5$ flat plumbings, $(c)$ a flat $5$-dipole surface of the link $4_1^2$ and
$(d)$ a graph diagram $D(K_{2,6})$ of a complete bipartite graph $K_{2,6}$ whose boundary is the link $4_1^2$, where the voltage assignments are all zero.} \label{412complete}
\end{figure}
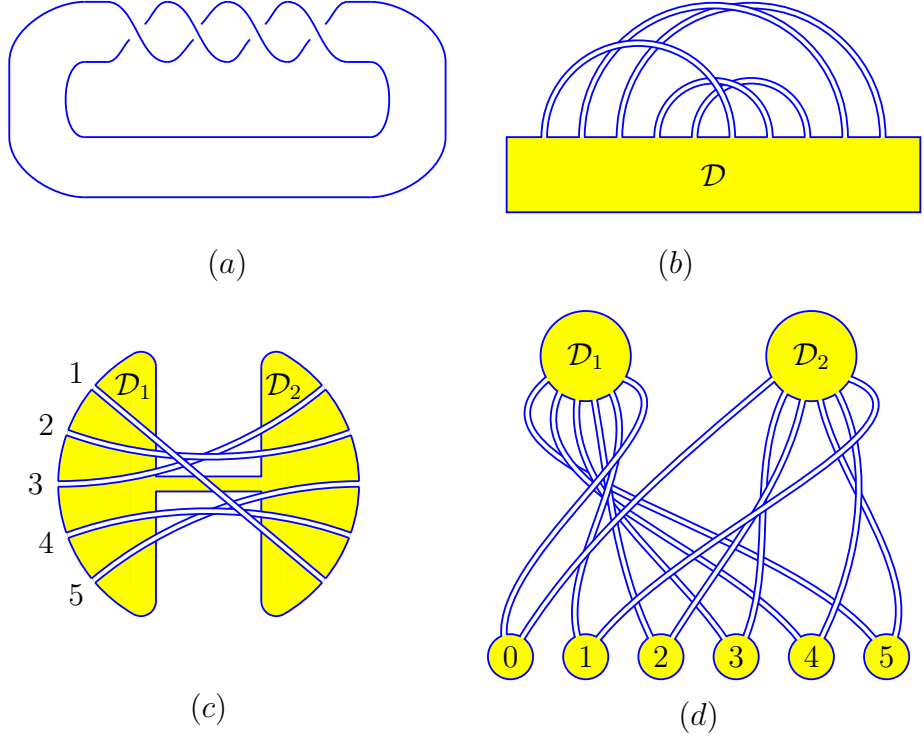

\vskip .5cm
\noindent{\tt Algorithm}
\begin{itemize}
\item $Step~1$. For a give link $L$, we find its braid representation $\beta$, the closed braid $\overline{\beta} =L$.
\item $Step~2$. Apply the method in~\cite{FHK:openbook} in order to obtain a flat plumbing basket surface $\F$ which is obtained from a disc by successively plumbing flat annuli.
\item $Step~3$. Apply the claim in Theorem~\ref{dipolemainthm} in order to find a flat plumbing basket surface $\F$ of $n$ annuli whose boundary is $L$, which
can be presented by two sets of $n$-tuples composed of $\{ 1, 2, \ldots, n\}$, the first $n$-tuples representing
how the bands in the flat dipole surface are connected and the second $n$-tuples representing the order of bands from the top to the bottom.
\item $Step~4$. Apply the move in Figure~\ref{dipolefig1}, and we obtain a flat $(n+1)$-dipole surface.
\item $Step~5$. Subdivide the bands by adding a disc in the middle of each band, as described in Corollary~\ref{dipolemaincor};
we obtain the desired graph diagram $D(K_{2,n+1})$ whose boundary is the given link $L$.
\end{itemize}

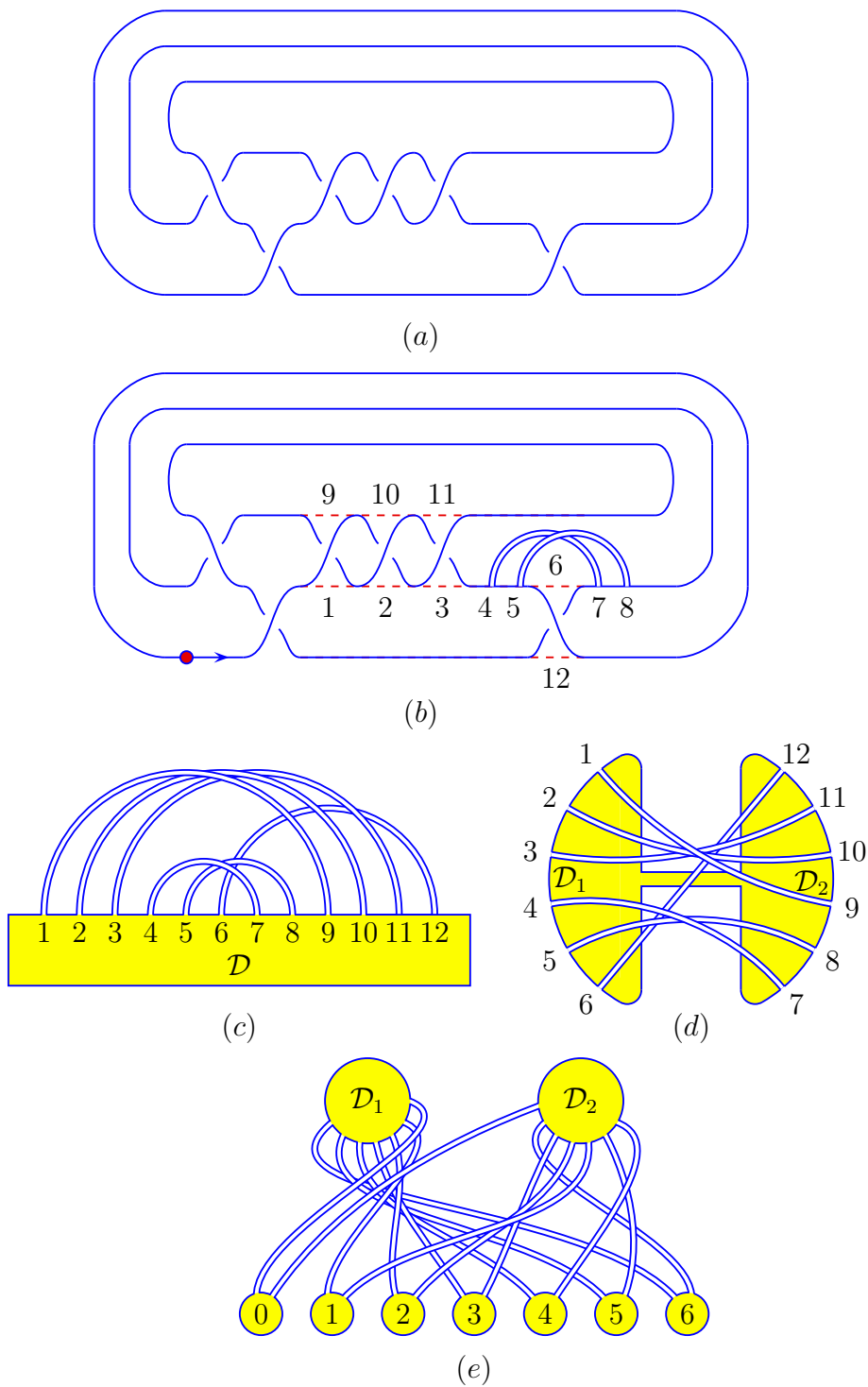
\begin{figure}
$$
\begin{pspicture}[shift=-1.2](0,-1.8)(9.5,3.1)
\psline(0,0)(0,2)
\pccurve[angleA=90,angleB=180](0,2)(1,3)
\psline(1,3)(8.2,3)
\pccurve[angleA=0,angleB=90](8.2,3)(9.2,2)
\psline(9.2,2)(9.2,0)
\pccurve[angleA=-90,angleB=0](9.2,0)(8.2,-1)
\psline(8.2,-1)(6.9,-1)
\pccurve[angleA=180,angleB=-45](6.9,-1)(6.6,-.6)
\pccurve[angleA=135,angleB=0](6.4,-.4)(6.1,0)
\psline(6.1,0)(5.3,0)
\pccurve[angleA=180,angleB=-45](5.3,0)(5,.4)
\pccurve[angleA=135,angleB=0](4.8,.6)(4.5,1)
\pccurve[angleA=180,angleB=0](4.5,1)(3.7,0)
\pccurve[angleA=180,angleB=-45](3.7,0)(3.4,.4)
\pccurve[angleA=135,angleB=0](3.2,.6)(2.9,1)
\psline(2.9,1)(2.1,1)
\pccurve[angleA=180,angleB=45](2.1,1)(1.8,.6)
\pccurve[angleA=-135,angleB=0](1.6,.4)(1.3,0)
\psline(1.3,0)(1,0)
\pccurve[angleA=180,angleB=-90](1,0)(.5,.5)
\psline(.5,.5)(.5,2)
\pccurve[angleA=90,angleB=180](.5,2)(1,2.5)
\psline(1,2.5)(8.2,2.5)
\pccurve[angleA=0,angleB=90](8.2,2.5)(8.7,2)
\psline(8.7,2)(8.7,.5)
\pccurve[angleA=-90,angleB=0](8.7,.5)(8.2,0)
\psline(8.2,0)(6.9,0)
\pccurve[angleA=180,angleB=0](6.9,0)(6.1,-1)
\psline(6.1,-1)(2.9,-1)
\pccurve[angleA=180,angleB=-45](2.9,-1)(2.6,-.6)
\pccurve[angleA=135,angleB=0](2.4,-.4)(2.1,0)
\pccurve[angleA=180,angleB=0](2.1,0)(1.3,1)
\pccurve[angleA=180,angleB=180](1.3,1)(1.3,2)
\psline(1.3,2)(7.9,2)
\pccurve[angleA=0,angleB=0](7.9,2)(7.9,1)
\psline(7.9,1)(5.3,1)
\pccurve[angleA=180,angleB=0](5.3,1)(4.5,0)
\pccurve[angleA=180,angleB=-45](4.5,0)(4.2,.4)
\pccurve[angleA=135,angleB=0](4,.6)(3.7,1)
\pccurve[angleA=180,angleB=0](3.7,1)(2.9,0)
\pccurve[angleA=180,angleB=0](2.9,0)(2.1,-1)
\psline(2.1,-1)(1,-1)
\pccurve[angleA=180,angleB=-90](1,-1)(0,0)
\rput(4.6,-1.6){{$(a)$}}
\end{pspicture}$$
$$
\begin{pspicture}[shift=-1.2](0,-2)(9.5,3.1)
\psline[linecolor=darkred, linestyle=dashed](2.9,1)(6.9,1)
\psline[linecolor=darkred, linestyle=dashed](2.9,0)(6.9,0)
\psline[linecolor=darkred, linestyle=dashed](2.9,-1)(6.9,-1)
\psline(0,0)(0,2)
\pccurve[angleA=90,angleB=180](0,2)(1,3)
\psline(1,3)(8.2,3)
\pccurve[angleA=0,angleB=90](8.2,3)(9.2,2)
\psline(9.2,2)(9.2,0)
\pccurve[angleA=-90,angleB=0](9.2,0)(8.2,-1)
\psline(8.2,-1)(6.9,-1)
\pccurve[angleA=180,angleB=0](6.9,-1)(6.1,0)
\psline(6.1,0)(5.3,0)
\pccurve[angleA=180,angleB=-45](5.3,0)(5,.4)
\pccurve[angleA=135,angleB=0](4.8,.6)(4.5,1)
\pccurve[angleA=180,angleB=0](4.5,1)(3.7,0)
\pccurve[angleA=180,angleB=-45](3.7,0)(3.4,.4)
\pccurve[angleA=135,angleB=0](3.2,.6)(2.9,1)
\psline(2.9,1)(2.1,1)
\pccurve[angleA=180,angleB=45](2.1,1)(1.8,.6)
\pccurve[angleA=-135,angleB=0](1.6,.4)(1.3,0)
\psline(1.3,0)(1,0)
\pccurve[angleA=180,angleB=-90](1,0)(.5,.5)
\psline(.5,.5)(.5,2)
\pccurve[angleA=90,angleB=180](.5,2)(1,2.5)
\psline(1,2.5)(8.2,2.5)
\pccurve[angleA=0,angleB=90](8.2,2.5)(8.7,2)
\psline(8.7,2)(8.7,.5)
\pccurve[angleA=-90,angleB=0](8.7,.5)(8.2,0)
\psline(8.2,0)(6.9,0)
\pccurve[angleA=180,angleB=45](6.9,0)(6.6,-.4)
\pccurve[angleA=-135,angleB=0](6.4,-.6)(6.1,-1)
\psline(6.1,-1)(2.9,-1)
\pccurve[angleA=180,angleB=-45](2.9,-1)(2.6,-.6)
\pccurve[angleA=135,angleB=0](2.4,-.4)(2.1,0)
\pccurve[angleA=180,angleB=0](2.1,0)(1.3,1)
\pccurve[angleA=180,angleB=180](1.3,1)(1.3,2)
\psline(1.3,2)(7.9,2)
\pccurve[angleA=0,angleB=0](7.9,2)(7.9,1)
\psline(7.9,1)(5.3,1)
\pccurve[angleA=180,angleB=0](5.3,1)(4.5,0)
\pccurve[angleA=180,angleB=-45](4.5,0)(4.2,.4)
\pccurve[angleA=135,angleB=0](4,.6)(3.7,1)
\pccurve[angleA=180,angleB=0](3.7,1)(2.9,0)
\pccurve[angleA=180,angleB=0](2.9,0)(2.1,-1)
\psline(2.1,-1)(1,-1)
\pccurve[angleA=180,angleB=-90](1,-1)(0,0)
\pscircle[linecolor=blue, fillcolor=darkred, fillstyle=solid](1.3,-1){.1}
\psline[arrowscale=1.5]{->}(1.8,-1)(1.9,-1)
\psarc[doubleline=true](6.35,0){.75}{0}{180}
\psarc[doubleline=true](6.75,0){.75}{0}{180}
\rput(3.3,-.3){{$1$}}
\rput(4.1,-.3){{$2$}}
\rput(4.9,-.3){{$3$}}
\rput(5.5,-.3){{$4$}}
\rput(5.9,-.3){{$5$}}
\rput(6.5,.3){{$6$}}
\rput(7.1,-.3){{$7$}}
\rput(7.5,-.3){{$8$}}
\rput(3.3,1.3){{$9$}}
\rput(4.1,1.3){{$10$}}
\rput(4.9,1.3){{$11$}}
\rput(6.5,-1.3){{$12$}}
\rput(4.6,-1.8){{$(b)$}}
\end{pspicture}$$
$$
\begin{pspicture}[shift=-.8](-.7,-2.4)(6.2,2)
\psarc[doubleline=true](4,-.5){1.5}{-5}{185}
\psarc[doubleline=true](2.75,-.5){.75}{-5}{185}
\psarc[doubleline=true](2.25,-.5){.75}{-5}{185}
\psarc[doubleline=true](3,-.5){2}{-5}{185}
\psarc[doubleline=true](2.5,-.5){2}{-5}{185}
\psarc[doubleline=true](2,-.5){2}{-5}{185}
\psframe[linecolor=lightgray,fillstyle=solid,fillcolor=lightgray](-.5,-1.5)(6,-.5)
\psline(-.03,-.5)(-.5,-.5)(-.5,-1.5)(6,-1.5)(6,-.5)(4.53,-.5)
\psline(.03,-.5)(.47,-.5) \psline(.53,-.5)(.97,-.5) \psline(1.03,-.5)(1.47,-.5)
\psline(1.53,-.5)(1.97,-.5) \psline(2.03,-.5)(2.47,-.5) \psline(2.53,-.5)(2.97,-.5)
\psline(3.03,-.5)(3.47,-.5) \psline(3.53,-.5)(3.97,-.5) \psline(4.03,-.5)(4.47,-.5)
\psline(4.53,-.5)(4.97,-.5) \psline(5.03,-.5)(5.47,-.5)
\rput(0,-.75){{$1$}}
\rput(.5,-.75){{$2$}}
\rput(1,-.75){{$3$}}
\rput(1.5,-.75){{$4$}}
\rput(2,-.75){{$5$}}
\rput(2.5,-.75){{$6$}}
\rput(3,-.75){{$7$}}
\rput(3.5,-.75){{$8$}}
\rput(4,-.75){{$9$}}
\rput(4.5,-.75){{$10$}}
\rput(5,-.75){{$11$}}
\rput(5.5,-.75){{$12$}}
\rput(2.75,-1.2){{$\mathcal{D}$}}
\rput(2.75,-2.1){{$(c)$}}
\end{pspicture}
\quad
\begin{pspicture}[shift=-.8](-2.5,-2.4)(2.5,2)
\pscircle[linecolor=lightgray, fillcolor=lightgray, fillstyle=solid](0,0){2}
\psframe[linecolor=white, fillcolor=white, fillstyle=solid](-1,-2.1)(1,2.1)
\pscircle[linecolor=lightgray, fillcolor=lightgray, fillstyle=solid](1.8;120){.2}
\pscircle[linecolor=lightgray, fillcolor=lightgray, fillstyle=solid](1.8;60){.2}
\pscircle[linecolor=lightgray, fillcolor=lightgray, fillstyle=solid](1.8;240){.2}
\pscircle[linecolor=lightgray, fillcolor=lightgray, fillstyle=solid](1.8;-60){.2}
\psframe[linecolor=lightgray, fillcolor=lightgray, fillstyle=solid](-1,-1.5588)(-.7,1.5588)
\psframe[linecolor=lightgray, fillcolor=lightgray, fillstyle=solid](1,-1.5588)(.7,1.5588)
\psframe[linecolor=lightgray, fillcolor=lightgray, fillstyle=solid](-.7,-.1)(.7,.1)
\psline(-.7,1.5588)(-.7,.1)(.7,.1)(.7,1.5588) \psline(.7,-1.5588)(.7,-.1)(-.7,-.1)(-.7,-1.5588)
\pccurve[doubleline=true, angleA=50,angleB=-130](2;230)(2;50)
\pccurve[doubleline=true, angleA=30,angleB=150](2;210)(2;-30)
\pccurve[doubleline=true, angleA=10,angleB=130](2;190)(2;-50)
\pccurve[doubleline=true, angleA=-10,angleB=-150](2;170)(2;30)
\pccurve[doubleline=true, angleA=-30,angleB=-170](2;150)(2;10)
\pccurve[doubleline=true, angleA=-50,angleB=170](2;130)(2;-10)
\psarc(0,0){2}{11}{29} \psarc(0,0){2}{31}{49}
\psarc(0,0){2}{51}{60} \psarc(0,0){2}{120}{129}
\psarc(0,0){2}{131}{149} \psarc(0,0){2}{151}{169}
\psarc(0,0){2}{171}{189} \psarc(0,0){2}{191}{209}
\psarc(0,0){2}{211}{229} \psarc(0,0){2}{231}{240}
\psarc(0,0){2}{300}{309} \psarc(0,0){2}{351}{9}
\psarc(0,0){2}{311}{329} \psarc(0,0){2}{331}{349}
\psarc(1.8;60){.2}{60}{180}
\psarc(1.8;120){.2}{0}{120}
\psarc(1.8;240){.2}{240}{0}
\psarc(1.8;300){.2}{180}{300}
\rput(-1.7,0){{$\mathcal{D}_1$}}
\rput(1.7,-.05){{$\mathcal{D}_2$}}
\rput(2.3;130){{$1$}}
\rput(2.3;150){{$2$}}
\rput(2.3;170){{$3$}}
\rput(2.3;190){{$4$}}
\rput(2.3;210){{$5$}}
\rput(2.3;230){{$6$}}
\rput(2.3;-50){{$7$}}
\rput(2.3;-30){{$8$}}
\rput(2.3;-10){{$9$}}
\rput(2.3;10){{$10$}}
\rput(2.3;30){{$11$}}
\rput(2.3;50){{$12$}}
\rput(0,-2.1){{$(d)$}}
\end{pspicture} $$
$$
\begin{pspicture}[shift=-1.2](-.8,-1)(5.6,3.6)
\pccurve[doubleline=true, angleA=120,angleB=-140](6,0)(1,2.7)
\pccurve[doubleline=true, angleA=60,angleB=-140](6,0)(4,2.7)
\pccurve[doubleline=true, angleA=120,angleB=-120](5,0)(1.2,2.6)
\pccurve[doubleline=true, angleA=60,angleB=-60](5,0)(4.8,2.6)
\pccurve[doubleline=true, angleA=120,angleB=-100](4,0)(1.4,2.5)
\pccurve[doubleline=true, angleA=60,angleB=-20](4,0)(5,2.7)
\pccurve[doubleline=true, angleA=120,angleB=-80](3,0)(1.6,2.5)
\pccurve[doubleline=true, angleA=60,angleB=-120](3,0)(4.2,2.6)
\pccurve[doubleline=true, angleA=120,angleB=-60](2,0)(1.8,2.6)
\pccurve[doubleline=true, angleA=60,angleB=-100](2,0)(4.4,2.5)
\pccurve[doubleline=true, angleA=120,angleB=-20](1,0)(2,2.7)
\pccurve[doubleline=true, angleA=60,angleB=-80](1,0)(4.6,2.5)
\pccurve[doubleline=true, angleA=120,angleB=-20](0,0)(2.08,3.0)
\pccurve[doubleline=true, angleA=60,angleB=-160](0,0)(4.1,3.0)
\pscircle[linecolor=lightgray, fillcolor=lightgray, fillstyle=solid](0,0){.3}
\pscircle[linecolor=lightgray, fillcolor=lightgray, fillstyle=solid](1,0){.3}
\pscircle[linecolor=lightgray, fillcolor=lightgray, fillstyle=solid](2,0){.3}
\pscircle[linecolor=lightgray, fillcolor=lightgray, fillstyle=solid](3,0){.3}
\pscircle[linecolor=lightgray, fillcolor=lightgray, fillstyle=solid](4,0){.3}
\pscircle[linecolor=lightgray, fillcolor=lightgray, fillstyle=solid](5,0){.3}
\pscircle[linecolor=lightgray, fillcolor=lightgray, fillstyle=solid](6,0){.3}
\pscircle[linecolor=lightgray, fillcolor=lightgray, fillstyle=solid](1.5,3){.6}
\pscircle[linecolor=lightgray, fillcolor=lightgray, fillstyle=solid](4.5,3){.6}
\psarc(0,0){.3}{108}{53} \psarc(0,0){.3}{63}{96}
\psarc(1,0){.3}{108}{48} \psarc(1,0){.3}{58}{96}
\psarc(2,0){.3}{117}{50} \psarc(2,0){.3}{61}{106}
\psarc(3,0){.3}{127}{57} \psarc(3,0){.3}{67}{117}
\psarc(4,0){.3}{132}{52} \psarc(4,0){.3}{63}{121}
\psarc(5,0){.3}{133}{60.5} \psarc(5,0){.3}{71}{122}
\psarc(6,0){.3}{133}{72} \psarc(6,0){.3}{85}{122}
\psarc(1.5,3){.6}{-28}{-3.5}
\psarc(1.5,3){.6}{2}{209} \psarc(1.5,3){.6}{214}{232}
\psarc(1.5,3){.6}{237.5}{256.5} \psarc(1.5,3){.6}{262}{278.5}
\psarc(1.5,3){.6}{284}{302.5} \psarc(1.5,3){.6}{308}{327}
\psarc(4.5,3){.6}{-28}{184}
\psarc(4.5,3){.6}{190}{209} \psarc(4.5,3){.6}{214}{232}
\psarc(4.5,3){.6}{237.5}{256.5} \psarc(4.5,3){.6}{262}{278.5}
\psarc(4.5,3){.6}{284}{302.5} \psarc(4.5,3){.6}{308}{327}
\rput(3,-.8){{$(e)$}}
\rput(1.5,3){{$\mathcal{D}_1$}}
\rput(4.5,3){{$\mathcal{D}_2$}}
\rput(0,0){{$0$}}
\rput(1,0){{$1$}}
\rput(2,0){{$2$}}
\rput(3,0){{$3$}}
\rput(4,0){{$4$}}
\rput(5,0){{$5$}}
\rput(6,0){{$6$}}
\end{pspicture}
$$
\caption{$(a)$ The knot $5_2$ as a closed braid, $(b)$ Seifert surface of $5_2$ in order to apply the algorithm in~\cite{FHK:openbook},
$(c)$ a flat plumbing basket surface of $5_2$, $(d)$ a flat $6$-dipole surface of $5_2$ and
$(e)$ a graph diagram $D(K_{2,7})$ whose boundary is $5_2$ where the voltage assignments are all zero.} \label{52complete}
\end{figure}

Now, we apply the algorithm for the links $5_2$ and $4_1^2$ in the following examples.

\begin{exa} \label{exa1}
A graph diagram $D(K_{2,6})$ of a complete bipartite graph $K_{2,6}$ whose boundary is the link $4_1^2$ where
the voltage assignments are all zero as illustrated in Figure~\ref{412complete} $(d)$.
\begin{proof}
The link $4_1^2$ is a closed braid $(\sigma_1)^4$ of two strings, as drawn in
Figure~\ref{412complete} $(a)$. In ~\cite{FHK:openbook}, it was shown that
$4_1^2$ has a flat plumbing surface with $3$ flat plumbings and a flat plumbing basket surface with $5$ flat plumbings as illustrated in
Figure~\ref{412complete} $(b)$. This one already satisfies the hypothesis of a flat dipole surface; thus, we have $<(1, 4, 5, 2, 3) | (1,2,3,4,5)>$
to present it as described in $Step~3$. We shrink the disc $\D$ into the union of the two discs to have a flat $6$-dipole surface, as depicted as
the shaded region in Figure~\ref{412complete} $(c)$. By following $Step~4$, we obtain a graph diagram $D(K_{2,6})$ of
a complete bipartite graph $K_{2,6}$ whose boundary is the link $4_1^2$ where the voltage assignments are all zero as shown in Figure~\ref{412complete} $(d)$.
\end{proof}
\end{exa}

\begin{exa} \label{exa2}
A graph diagram $D(K_{2,7})$ of a complete bipartite graph $K_{2,7}$ whose boundary is the knot $5_2$ where
the voltage assignments are all zero as illustrated in Figure~\ref{52complete}.
\begin{proof}
The link $5_2$ is a closed braid $\sigma_2^{-1}\sigma_1(\sigma_2)^3 \sigma_1$ of three strings as drawn in
Figure~\ref{52complete} $(a)$. \cite{FHK:openbook} $(a)$ provided an algorithm to find a link's
flat plumbing surface as depicted in Figure~\ref{52complete}. Let us explain the algorithm. First we choose the disc $\D$ which is the union of Seifert discs connected by
two half twisted bands which are represented by $\sigma_2^{-1}\sigma_1$ as drawn in the figure by the dashed red line. For possibility of flat plumbing, we add two extra annuli.
The numbers $1, 2, \ldots, 12$ in the figure were chosen by reading the order of flat plumbing from the point along the direction of the arrow in the figure.
The resulting flat plumbing basket surface with $6$ flat plumbings are given in
Figure~\ref{52complete} $(c)$. This one already satisfies the hypothesis of a flat dipole surface which can be read as $<(4, 5, 1, 2, 3, 6) | (1,2,3,4,5,6)>$
in order to present it as described in $Step~3$. We shrink the disc $\D$ to the union of two discs $\D_1$ and $D_2$ to have a flat $6$-dipole surface as depicted as
the shaded region in Figure~\ref{52complete} $(d)$. By following $Step~4$, we obtain a graph diagram $D(K_{2,6})$ of
a complete bipartite graph $K_{2,7}$ whose boundary is the knot $5_2$ where the voltage assignments are all zero as shown in Figure~\ref{52complete} $(e)$.
\end{proof}
\end{exa}

\section*{Acknowledgments}
The author would like to thank Younghae Do, Hunki Baek and Myoungsoo Seo for their
helpful discussions. The \TeX\, macro package
PSTricks~\cite{PSTricks} was essential for typesetting the equations
and figures. This work was supported by Kyonggi University Research Grant 2011.

\end{document}